\documentclass[reqno]{amsart}
\usepackage{amsmath,amssymb}
\newcommand{\al}{\ensuremath{\alpha}}
\newcommand{\be}{\ensuremath{\beta}}

\newcommand{\de}{\ensuremath{\delta}}
\newcommand{\ep}{\ensuremath{\varepsilon}}

\newcommand{\la}{\ensuremath{\lambda}}
\newcommand{\si}{\ensuremath{\sigma}}

\newcommand{\Ga}{\ensuremath{\Gamma}}

\newcommand{\CC}{\ensuremath{\mathbb{C}}}

\newcommand{\RR}{\ensuremath{\mathbb{R}}}
\newcommand{\TT}{\ensuremath{\mathbb{T}}}
\newcommand{\ZZ}{\ensuremath{\mathbb{Z}}}

\newcommand{\FSH}{\ensuremath{\mathcal H}}
\newcommand{\FSI}{\ensuremath{\mathcal I}}

\newcommand{\FSW}{\ensuremath{\mathcal W}}

\newcommand{\FSZ}{\ensuremath{\mathcal Z}}

\newcommand{\goa}{\ensuremath{\mathfrak{a}}}
\newcommand{\goh}{\ensuremath{\mathfrak{h}}}
\newcommand{\gog}{\ensuremath{\mathfrak{g}}}

\newcommand{\gok}{\ensuremath{\mathfrak{k}}}
\newcommand{\gol}{\ensuremath{\mathfrak{l}}}

\newcommand{\gop}{\ensuremath{\mathfrak{p}}}

\newcommand{\gou}{\ensuremath{\mathfrak{u}}}

\newcommand{\goS}{\ensuremath{\mathfrak{S}}}

\DeclareMathOperator{\id}{id}

\DeclareMathOperator{\End}{End}
\DeclareMathOperator{\Hom}{Hom}\DeclareMathOperator{\sgn}{sgn_q}
\DeclareMathOperator{\Mat}{Mat}
\DeclareMathOperator{\pr}{pr}

\newcommand{\Aq}{\ensuremath{A_q}}
\newcommand{\AG}{\ensuremath{A_q(G)}}
\newcommand{\AK}{\ensuremath{A_q(K)}}
\newcommand{\AU}{\ensuremath{A_q(U)}}
\newcommand{\AUK}{\ensuremath{A_q(U/K)}}
\newcommand{\Aksi}{\ensuremath{A_q(\gok^\si\backslash U)}}
\newcommand{\Uq}{\ensuremath{U_q({\gog})}}
\newcommand{\detq}{\ensuremath{\det\nolimits_q}}
\newcommand{\ksi}{\ensuremath{\gok^\si}}
\newcommand{\ktau}{\ensuremath{\gok^\tau}}
\newcommand{\Hst}{\ensuremath{\FSH^{\si,\tau}}}
\newcommand{\Hinft}{\ensuremath{\FSH^{\infty,\tau}}}
\newcommand{\AT}{\ensuremath{A(\TT)}}
\newcommand{\Uh}{\ensuremath{U_q(\goh)}}
\newcommand{\ten}{\ensuremath{\otimes}}
\newcommand{\vek}[1]{\ensuremath{\mathbf{#1}}}
\newcommand{\phist}{\ensuremath{\varphi_r^{\si,\tau}}}

\numberwithin{equation}{section}
\theoremstyle{plain}
\newtheorem{thm}{Theorem}[section]
\newtheorem{cor}[thm]{Corollary}
\newtheorem{lem}[thm]{Lemma}
\newtheorem{prop}[thm]{Proposition}
\theoremstyle{remark}
\newtheorem{rem}[thm]{Remark}
\begin{document}
\title[$BC$ type orthogonal polynomials on quantum Grassmannians]{Some
limit transitions between $BC$ type \\ orthogonal polynomials 
interpreted on\\ quantum complex Grassmannians}
\author{Mathijs S. Dijkhuizen}
\author{Jasper V. Stokman}
\address{Department of Mathematics, Faculty of Science,
   Kobe University, Rokko, Kobe 657, Japan}
\email{msdz@math.kobe-u.ac.jp}
\address{Korteweg-de Vries institute for mathematics, University of Amsterdam, 
Plantage Muidergracht 24, 1018 TV Amsterdam, The Netherlands}
\email{jasper@wins.uva.nl}
\date{}
\begin{abstract}
The quantum complex Grassmannian $U_q/K_q$ of rank $l$ is
the quotient of the quantum unitary group $U_q=U_q(n)$ by the 
quantum subgroup $K_q=U_q(n-l)\times U_q(l)$. 
We show that $(U_q,K_q)$ is a quantum Gelfand pair and we 
express the zonal spherical functions, 
i.e.\ $K_q$-biinvariant matrix coefficients of finite-dimensional irreducible
representations of $U_q$, as multivariable little $q$-Jacobi
polynomials depending on one discrete parameter. Another type of
biinvariant matrix coefficients is identified as multivariable
big $q$-Jacobi polynomials. The proof is based on earlier results
by Noumi, Sugitani and the first author 
relating Koornwinder polynomials to a one-parameter
family of quantum complex Grassmannians, and certain limit
transitions from Koornwinder polynomials to multivariable big and little $q$-Jacobi
polynomials studied by Koornwinder and the second author.
\end{abstract}
\maketitle
%%%%%%%%%%%%%%%%%%%%%%%%%%%%%%%%%%%%%%%%%%%%%%%%%%%%%%%%%%%%%%%%%%%%%%%%%%%%%
%%
%%
%%  SECTION 1  INTRODUCTION
%%
%%
%%%%%%%%%%%%%%%%%%%%%%%%%%%%%%%%%%%%%%%%%%%%%%%%%%%%%%%%%%%%%%%%%%%%%%%%%%%%
\section{Introduction}
\label{section:intro}
The first connection between $q$-special functions and quantum groups
was revealed in the late 1980's by the interpretation of 
little $q$-Jacobi polynomials as matrix coefficients of irreducible 
representations of the quantum $SU(2)$ group
(cf.\  \cite{vs:su2}, \cite{masuda:su2}, \cite{koor:su2}). In the past decade,
many other connections between representation theory of quantum groups 
and the theory of $q$-special functions have been discovered.

For instance, Noumi and Mimachi \cite{NM} showed that
the zonal spherical functions on  Podle\`s's one-parameter
family of quantum 2-spheres can be identified with big $q$-Jacobi polynomials. 
In \cite{koornwinder} and \cite{koor:awsu2}, Koornwinder  
generalized these results by replacing the notion of invariance 
under a quantum subgroup by the notion of invariance under a twisted primitive 
element in the quantized universal enveloping algebra. 
In particular, this infinitesimal approach allowed Koornwinder to identify
the zonal spherical functions on Podle\`s spheres
as a two-parameter family of Askey-Wilson polynomials.
Noumi's and Mimachi's results \cite{NM} could then be 
reobtained by sending one parameter to infinity. If both parameters
are sent to infinity, one obtains the interpretation of little $q$-Jacobi polynomials
as zonal spherical functions on the ``standard'' $2$-sphere, i.e.\ the quantum
$2$-sphere which is realized as the quotient of the quantum $SU(2)$ group by
the standard maximal torus. On the level of $q$-special functions, 
these limits correspond to certain explicit limit transitions 
from Askey-Wilson polynomials to big respectively 
little $q$-Jacobi polynomials.

Analogous statements are valid for complex projective space, see 
Noumi, Yamada and Mimachi \cite{nym:gl} for the little $q$-Jacobi case,
and Dijkhuizen and Noumi \cite{dijk-nou:proj} for the general case.
In this paper we generalize these results to the higher rank setting
by interpreting certain subfamilies of Koornwinder's multivariable 
analogues of the Askey-Wilson polynomials and certain subfamilies of 
the multivariable big and little $q$-Jacobi polynomials (cf.\  \cite{stok:jac})
as zonal spherical functions on quantum analogues of the complex Grassmannian
\[
U/K:=U(n)/(U(n-l)\times U(l)),\qquad (l\leq \lbrack n/2\rbrack),
\]
where $U(n)$ is the group of $n\times n$ unitary matrices.

Koornwinder's infinitesimal approach
to harmonic analysis on quantized symmetric spaces was for the first time
successfully generalized to higher rank cases by Noumi \cite{nou:macdonald}. 
The quantized symmetric spaces were now defined using invariance under 
certain two-sided coideals in the quantized universal enveloping algebra. 
So far, this method has been successfully applied to all compact symmetric spaces of
classical type, e.g. \cite{nou-sug:toyonaka}, \cite{NS1}, \cite{nds:gras}, \cite{S},
\cite{dijk-nou:proj}. The related zonal spherical functions 
can all be identified with Koornwinder polynomials or Macdonald polynomials.

In particular, Noumi, Sugitani and the 
first author \cite{nds:gras} introduced a one-parameter family 
of quantum analogues of the complex Grassmannian.
They announced that the spherical functions associated with these quantized symmetric
spaces can be expressed as a two-parameter subfamily of the Koornwinder polynomials.

In this paper we extend these results to the quantum subgroup case, 
i.e. we determine the zonal spherical functions associated with the 
quantum analogue of the complex Grassmannian (denoted by $U_q/K_q$)
which is defined as the quotient of the quantum
unitary group $U_q(n)$ by the obvious quantum subgroup $K_q$
corresponding to $K=U(n-l)\times U(l)$. 
We show that the quantum space $U_q/K_q$
can be formally obtained from the one-parameter family of quantum Grassmannians
defined in \cite{nds:gras} by sending the parameter to infinity.
We also show that this limit transition on quantum Grassmannians 
is compatible with the limit transitions 
from Koornwinder polynomials to multivariable big and little $q$-Jacobi
polynomials, which were previously studied by Koornwinder and the second author
\cite{stok-koor:limit}. 

In order to give a rigorous meaning to the above-mentioned
formal limit between quantum Grassmannians 
we use the recent result of the second author \cite{stok:residue} that 
the limits from Koornwinder polynomials to multivariable big and little
$q$-Jacobi polynomials can be taken on the level of the orthogonality measures
(in a suitable weak sense). Indeed, this result allows us to rigorously compute  
the limit of the quantum Schur orthogonality relations 
for the zonal spherical functions. This in turn leads to a
complete description of the harmonic analysis on the quantum Grassmannian
$U_q/K_q$. In particular, we obtain an identification of 
the zonal spherical functions on the quantum complex Grassmannian $U_q/K_q$
with multivariable big respectively little $q$-Jacobi polynomials. 

We remark that in the rank one
case, there are several alternative methods 
for determining the zonal spherical functions in the quantum subgroup case
(see e.g. \cite{nym:gl}, \cite{VS2}).
However, the relatively indirect method using limit arguments
seems to be the only method which admits a direct generalization 
to the higher rank cases of the complex Grassmannian. 

There are strong indications that such limit arguments can also be applied
for other quantum compact symmetric spaces. 
Indeed, an important prerequisite for applying 
such limit arguments is the existence of a suitable one-parameter family 
of quantizations of the compact symmetric space under consideration.
The occurrence of such a one-parameter phenomena on the quantum level is directly
related to the existence of a one-parameter family of covariant
Poisson brackets on the underlying symmetric space. 
The existence of such a one-parameter family
of covariant Poisson brackets in the case of a Hermitian symmetric space was
established by  Khoroshkin \& Radul \& Rubtsov \cite{krr:herm} 
and Donin \& Gurevich \cite{donin-gur:Rmatrix}
(see also \cite{dijk:goslar} for more information and references). 

The paper is organized as follows.
In section
\ref{section:classical} we briefly recall some results on the
structure of the (classical) complex Grassmannian $U/K$ and 
on the nature of its zonal spherical functions.
In section \ref{section:polynomials} we recall the definition and the main
properties of the Koornwinder polynomials and the multivariable big and little 
$q$-Jacobi polynomials. In section \ref{section:preliminaries}
we collect some facts about the quantum unitary group that will
be heavily used in later sections.  In section \ref{section:spherical}
we determine the spherical dominant weights for the quantum complex Grassmannian
$U_q/K_q$. In section \ref{section:nds}
we recall the one-parameter family of quantum Grassmannians and determine
the corresponding spherical weights. In section 7 we express the zonal spherical
functions on the one-parameter family of quantum Grassmannians as
a subfamily of the Koornwinder polynomials. Many of the 
results in section \ref{section:nds} and section 7, 
previously announced in \cite{nds:gras}, are proved here 
in full detail using so-called  ``principal term'' type of arguments, 
cf.\  \cite{NS1}, \cite{S}. Finally, in section 8 we 
study the limit from the one-parameter family of quantum Grassmannians 
to $U_q/K_q$, which in particular leads to the interpretation of 
multivariable big and little $q$-Jacobi polynomials as 
zonal spherical functions on the quantum Grassmannian $U_q/K_q$.

Throughout this paper we will use the convention that vector spaces
are defined over the complex numbers and that algebras 
have a unit element.

%%%%%%%%%%%%%%%%%%%%%%%%%%%%%%%%%%%%%%%%%%%%%%%%%%%%%%%%%%%%%%%%%%%%%%%%%%%%%
%%
%%  SECTION 2  THE CLASSICAL COMPLEX GRASSMANNIAN 
%%
%%
%%%%%%%%%%%%%%%%%%%%%%%%%%%%%%%%%%%%%%%%%%%%%%%%%%%%%%%%%%%%%%%%%%%%%%%%%%%%%
\section{The classical complex Grassmannian}
\label{section:classical}
Two general references for the contents of this section are
Helgason \cite{hel:gga}, and Heckman and Schlichtkrull \cite{hs:harmonic}.

Throughout this paper, $n\geq 2$ and $1\leq l\leq [\frac{n}{2}]$ are
fixed integers. Let $G:=GL(n,\CC)$ denote the general linear group with Lie
algebra $\gog=\gog\gol(n,\CC)$, and $U:=U(n)$ the unitary group with
Lie algebra $\gou$.
Let $\TT\subset U$ denote the maximal torus consisting of diagonal
matrices in $U$. Write $\goh\subset \gog$ for the corresponding
Cartan subalgebra. Let $e_{ij}$ ($1\leq i,j\leq n$) denote the
standard matrix units. The matrices $h_i:= e_{ii}$ ($1\leq i\leq n$)
form a basis of $\goh$. Write $\tilde{\ep}_i\in \goh^\ast$ ($1\leq i\leq n$) 
for the corresponding dual basis vectors, and define a non-degenerate
symmetric bilinear form on ${\goh}^\ast$ by  $\langle
\tilde{\ep}_i,\tilde{\ep}_j\rangle=\delta_{ij}$.
The usual positive system $R^+$ in the root system $R:=R(\gog,\goh)$ 
consists of the vectors $\tilde{\ep}_i -\tilde{\ep}_j$ ($1\leq i<j\leq n$). 
Let $P=P_n:=\bigoplus_{1\leq i\leq n}\ZZ\tilde{\ep}_i$ denote the rational character
lattice of $G$ (equivalently, the lattice of analytically integral 
weights of $U$). 
Recall that the cone of dominant weights $P^+=P_n^+$ is given by
\begin{equation}
P_n^+:=\lbrace (\lambda_1,\ldots,\lambda_n)\in P\mid
\lambda_1\geq\lambda_2\geq\ldots\geq\lambda_n \rbrace.
\end{equation}
Denote by $\leq$ the (partial) dominance ordering on $P$. One has
$\mu\leq \lambda$ if and only if 
\begin{equation}\label{e:Adominance}
\sum_{i=1}^j\mu_i\leq\sum_{i=1}^j\lambda_i
\quad (1\leq j\leq n-1)\quad\text{and}\quad \sum_{i=1}^n \mu_i =
\sum_{i=1}^n \lambda_i.
\end{equation}

We write $K:=U(n-l)\times U(l)$ and $\gok:= \gog\gol(n-l,\CC) \oplus
\gog\gol(l,\CC)$ for the corresponding complexified Lie algebra $\gok$. 
$K$ is regarded as a subgroup of $U$ via the embedding
\begin{equation}\label{e:Kembedding}
U(n-l)\times U(l)\hookrightarrow U(n), \quad (A,B) \mapsto 
\begin{pmatrix}A & 0 \\ 0 &B \end{pmatrix}.
\end{equation}

The pair $(U,K)$ is symmetric. Indeed, 
the involutive Lie group automorphism $\theta: U\to U$
defined by $\theta(g):=JgJ$ with
\begin{equation}\label{e:Jdef}
J:=\sum_{1\leq k\leq n-l} e_{kk} - \sum_{1\leq k\leq l} e_{k'k'},\quad (k':=n+1-k)
\end{equation}
has fixed point group $K$. The differential of $\theta$ at the unit element $e\in U$,
extended $\CC$-linearly to a Lie algebra involution of $\gog$, will also be denoted
by $\theta$.
The $+1$ eigenspace of the 
involution $\theta\colon \gog\to\gog$ is exactly the Lie subalgebra $\gok$.
Writing $\gop$ for the $-1$ eigenspace of $\theta$ we have the eigenspace
decomposition $\gog=\gok\oplus \gop$.  

For certain purposes, it is more convenient to consider the involution
$\theta'\colon \gog\to\gog$ defined by $\theta'(X):= J'XJ'$ with 
\begin{equation}\label{e:Jprimedef}
J':=\sum_{l<k<l'} e_{kk} - \sum_{1\leq k\leq l} e_{kk'}
- \sum_{1\leq k\leq l} e_{k'k}.
\end{equation}
Since $J'$ is conjugate to $J$, the involution $\theta'$ is conjugate to 
$\theta$ by an inner automorphism of $\gog$. Observe that both $\theta$ and
$\theta'$ leave $\gou\subset \gog$ invariant. 

Let $\gog=\gok'\oplus\gop'$ be the
eigenspace decomposition of $\theta'$ into $+1$ and $-1$ eigenspaces.
The intersection $\goh\cap \gok'$
is spanned by the elements $h_i+h_{i'}$ ($1\leq i\leq l$) and
$h_i$ ($l< i < l'$), whereas the intersection $\goa:= \goh\cap \gop'$
is spanned by $h_i-h_{i'}$ ($1\leq i\leq l$) and is maximal abelian
in $\gop'$. The positive system of $R$ taken with respect to
the lexicographic ordering of $\goh_{\RR}:=\sum_{j=1}^n{\RR}h_j$
relative to the ordered basis
$h_1-h_{1'}, \ldots, h_l-h_{l'}$, $h_1+h_{1'}, \ldots, h_l+h_{l'}$,
$h_{l+1}, \ldots, h_{n-l}$ coincides with $R^+$.

Write $\ep'_i$ for the restriction of $\tilde{\ep}_i$ to $\goa$
($1\leq i\leq l$). The root system $R\subset \goh^\ast$ is mapped
under the natural projection $\goh^\ast \twoheadrightarrow \goa^\ast$ 
onto the restricted root system $\Sigma'=\Sigma'(\gog,\goa)$.
Choose the positive system in $\Sigma'$ with respect to the 
lexicographic ordering of $\goa_{\RR}:= \goh_\RR \cap \goa$ 
relative to the ordered basis $h_1-h_{1'}, \ldots, h_l-h_{l'}$ of $\goa_\RR$.
This ordering is compatible with the lexicographic ordering
of $\goh_\RR$ introduced above in the sense that $\lambda\in \goh^\ast_\RR$
is positive if its restriction to $\goa_\RR$ is strictly positive. The positive
root vectors in $\Sigma'$ are
\begin{equation*}
\ep'_i\quad (1\leq i\leq l), \quad
\ep'_i \pm \ep'_j \quad (1\leq i < j\leq l), \quad
2\ep'_i\quad (1\leq i\leq l),
\end{equation*}
the roots $\ep'_i$ ($1\leq i\leq l$) occurring only if $n\neq 2l$.
$\Sigma'$ is isomorphic with $BC_l$ if $n\neq 2l$ and
isomorphic with $C_l$ if $n=2l$.
The root multiplicities corresponding to the short, medium, and long
roots are
\begin{equation}\label{e:rootmult}
m_1 = 2(n-2l), \quad m_2 = 2\; (l>1), \quad m_3 = 1.
\end{equation}
For later purposes, it is convenient to rescale the root system $\Sigma'$ by a
factor $2$. So we set $\Sigma:=2\Sigma'\subset \goa^\ast$, $\ep_i := 2\ep'_i$
($1\leq i\leq l$). Then the corresponding weight lattice 
$P_\Sigma\subset \goa^\ast$ is
the $\ZZ$-span of the $\ep_i$ ($1\leq i\leq l$), 
and the set $P_\Sigma^+$ of
dominant weights $\mu =\sum_i \mu_i \ep_i$
(taken with respect to the lexicographic ordering on $\goa_\RR$ introduced above)
is characterized by the condition $\mu_1\geq \cdots \geq \mu_l \geq 0$.
The dominance ordering $\leq$ on $P_\Sigma$ is explicitly given by
\begin{equation}\label{e:BCdominance}
\mu\leq\lambda \Longleftrightarrow \sum_{i=1}^j\mu_i\leq\sum_{i=1}^j\lambda_i
\quad (1\leq j\leq l).
\end{equation}

Let $K'\subset U$ denote the connected subgroup corresponding to $\gok'$.
The symmetric pairs $(U,K)$ and $(U,K')$ are Gelfand pairs, i.e.\
every finite-dimensional irreducible representation of $U$ has at
most one $K$-fixed vector up to scalar multiples. According to
\cite[Chapter V, Theorem 4.1]{hel:gga}, a highest weight $\lambda\in P^+$
is $K'$-spherical, i.e.\ corresponds to a representation with a non-zero
$K'$-fixed vector, if and only if the restriction of $\lambda$ to
$\goh\cap \gok'$ is zero and the restriction of $\lambda$ to
$\goa$ lies in $P^+_\Sigma$. Hence we get the following result.
\begin{thm}\label{th:cl-gelfand}
The set $P^+_K\subset P^+$ of $K$-spherical dominant weights consists
of all dominant weights of the form
\begin{equation*}
\lambda:= (\lambda_1, \ldots, \lambda_l, 0,\ldots, 0, 
-\lambda_l,\ldots, -\lambda_1).
\end{equation*}
\end{thm}
Write $\lambda^\natural := (\lambda_1, \ldots, \lambda_l)$ for a dominant
weight $\lambda\in P^+_K$. The assignment $\lambda\mapsto \lambda^\natural$ 
defines a bijection of $P^+_K$ onto $P_\Sigma^+$.
Let $\varpi_r\in P_K^+$ be the spherical weight for which 
$\varpi_r^\natural=(1^r)$. Then $P_K^+=\bigoplus_{1\leq r\leq
l}{\mathbb{Z}}_+\varpi_r$. We will call $\{ \varpi_r\}_{r=1}^l$ the 
fundamental dominant spherical weights.

Let $A$ denote the algebra of polynomial functions on $U$,
$\AT$ the algebra of polynomial functions on the maximal torus $\TT$.
$\AT$ may be naturally identified with the algebra $\CC[z^{\pm 1}]=\CC[z_1^{\pm 1},
\ldots,z_n^{\pm 1}]$ in $n$ variables $z_i$ ($1\leq i\leq n$) in the following way. 
Observe that $\TT\simeq i\goh_{\mathbb{R}}/2\pi i \tilde{P}$ 
via the exponential mapping, 
where $\tilde{P}:=\oplus_{1\leq j \leq n}{\mathbb{Z}}h_j$. 
Then the coordinate functions $z_j$ can be defined by $z_j:=e^{\tilde{\ep}_j}$, 
where $e^{\tilde{\ep}_j}([X]):=e^{\tilde{\ep}_j(X)}$
for $[X]\in T$ with $X\in i\goh_{\mathbb{R}}$ a representative of $[X]$.
More explicitly, $z_j$ is given by
\begin{equation*}
z_j: \hbox{diag}(e^{i\theta_1},e^{i\theta_2},\ldots,e^{i\theta_n})\mapsto
e^{i\theta_j},\quad (\theta_k\in [0,2\pi))
\end{equation*}
where $\hbox{diag}(a_1,\ldots,a_n)$ is the diagonal matrix with $a_1,\ldots,a_n$
on the diagonal.

Let $\FSH\subset A$ denote the subalgebra of $K'$-biinvariant functions. 
One has the decomposition
\begin{equation}\label{e:cl-Hdecomp}
\FSH = \bigoplus_{\lambda\in P^+_K} \FSH(\lambda), \quad
\FSH(\lambda) := \FSH \cap W(\lambda),
\end{equation}
$W(\lambda)\subset A$ denoting the subspace spanned by the matrix coefficients
of the irreducible representation of highest weight $\lambda$. Each
of the subspaces $\FSH(\lambda)$ ($\lambda\in P^+_K$) is one-dimensional,
since $(U,K')$ is a Gelfand pair.
Any non-zero element $\varphi(\lambda)$ 
of $\FSH(\lambda)$ is called a zonal spherical function.

Set
\begin{equation}
\TT_l := \hbox{exp}(i\goa_{\mathbb{R}})/(\hbox{exp}(i\goa_{\mathbb{R}})\cap K')
\simeq i\goa_\RR/2\pi iQ^\vee_\Sigma,
\end{equation} 
with $Q^\vee_\Sigma\subset \goa_\RR$ the
coroot lattice of $\Sigma$. More explicitly, 
the coroot lattice $Q^\vee_\Sigma$ is the $\ZZ$-span
of the elements $\frac{1}{2}(h_i-h_{i'})$ ($1\leq i\leq l$),
$\hbox{exp}(i\goa_{\mathbb{R}})$ are the diagonal matrices
\[
\hbox{diag}(e^{i\theta_1},\ldots e^{i\theta_l},1,\ldots,1,e^{-i\theta_l},\ldots,
e^{-i\theta_1})\quad (\theta_j\in [0,2\pi)),
\]
and $\hbox{exp}(i\goa_{\mathbb{R}})\cap K'$ are the matrices in
$\hbox{exp}(i\goa_{\mathbb{R}})$ of order $2$.

Write $\log: \TT_l\rightarrow i\goa_{\mathbb{R}}$ for 
the multi-valued inverse of the exponential map $\exp: i\goa_{\mathbb{R}}\to
\TT_l$.  Similarly as for $\AT$, 
the algebra of polynomial functions on $\TT_l$ may be identified 
with the algebra $\CC[x^{\pm 1}]$ of Laurent polynomials in 
the $l$ variables $x_j$ ($1\leq j\leq l$), where the identification is
given by $x_j(t):=e^{\ep_j(\log(t))}$. In other words, the map $x_j$ 
is given by
\begin{equation*}
x_j: \hbox{diag}(e^{i\theta_1},\ldots e^{i\theta_l},1,\ldots,1,e^{-i\theta_l},\ldots,
e^{-i\theta_1})\mapsto e^{2i\theta_j}.
\end{equation*}

It follows that the algebra ${\mathbb{C}}[x^{\pm 1}]$ of polynomial 
functions on $\TT_l$ can be naturally embedded in 
the algebra ${\mathbb{C}}[z^{\pm 1}]$ of polynomial functions
on the maximal torus $\TT$ by the assignment
\begin{equation}\label{e:xdef}
x_1 = z_1z_n^{-1}, \quad x_2= z_2z_{n-1}^{-1},
\quad \ldots, \quad x_l = z_l z_{n+1-l}^{-1}.
\end{equation}
Let $\goS:=\goS_l$ denote the permutation group on $l$ letters,
$\FSW:=\FSW_l=\ZZ_2^l \rtimes \goS_l$ the Weyl group of $\Sigma$. The natural action
of $\FSW$ on $\goa_\RR$ descends to $\TT_l$. Hence $\FSW$ acts
naturally on the algebra $\CC[x^{\pm 1}]$.
Write $\CC[x^{\pm 1}]^\FSW$ for the subalgebra of $\FSW$-invariant 
Laurent polynomials. By Chevalley's restriction theorem and the above-mentioned
natural embedding of $\CC[x^{\pm 1}]$ into $\CC[z^{\pm 1}]$, we have the following
theorem.
\begin{thm}\label{th:cl-Winvar}
Restriction to $\TT$ induces an isomorphism
of $\FSH$ onto the algebra $\CC[x^{\pm 1}]^\FSW$ of $\FSW$-invariant
Laurent polynomials in the variables $x_i$ \textup{(}$1\leq i\leq l$\textup{)}.
\end{thm}

By Theorem \ref{th:cl-Winvar} the direct sum decomposition 
\eqref{e:cl-Hdecomp} of $\FSH$ gives rise to a unique (up to rescaling) 
linear basis of $\CC[x^{\pm 1}]^{\FSW}$.
This linear basis can be expressed in terms of
$BC$ type Heckman-Opdam polynomials,  which we will now define.

Let $V_{HO}$ denote the set of triples $k:=(k_1,k_2,k_3)$ of real numbers
such that $k_1+k_3> -\frac{1}{2}$, $k_2>0$, $k_3> -\frac{1}{2}$. Define
an inner product 
$\langle \cdot, \cdot \rangle_{HO} = \langle \cdot,\cdot\rangle_{HO}^k$
on $\CC[x^{\pm 1}]^\FSW$ by
\begin{equation}\label{e:HOorth}
\langle P,Q\rangle_{HO} := \int_{\TT_l} P(t) \overline{Q(t)} \Delta_{HO}(t;k) dt.
\end{equation}
Here $dt$ denotes the normalized Haar measure on the 
torus $\TT_l$. The continuous positive weight function
$t\mapsto\Delta_{HO}(t;k)$ on $\TT_l$ is defined by
\begin{equation}
\Delta_{HO}(t;k) := \prod_{\al\in\Sigma} \left (e^{\frac{1}{2}i\langle \al,
\log(t)\rangle} -
e^{-\frac{1}{2}i\langle \al, \log(t)\rangle}\right )^{k_\al}.
\end{equation}
The multiplicity parameters
$k_\al$ are by definition  equal to $k_i$ ($i=1,2,3$), 
depending on whether $\al\in \Sigma$ is a short, medium, or long root.
If $l=1$ there is no dependence on $k_2$.

Recall that the usual orbit sums
\begin{equation}
\tilde{m}_{\lambda}(x):=\sum_{\mu\in {\FSW}\lambda}x^{\mu}
\quad (\lambda\in P_\Sigma^+)
\end{equation}
form a linear basis for the algebra $\CC[x^{\pm 1}]^\FSW$.
The Heckman-Opdam hypergeometric orthogonal polynomials
associated with the root system $\Sigma=BC_l$ (cf.\ \cite[Part I]{hs:harmonic})
are the uniquely determined family 
$\lbrace P^{HO}_\lambda \mid \lambda\in P^+_\Sigma\rbrace$ 
such that
\begin{equation}\label{e:defpolHO}
(i)\; P_\lambda^{HO} = \tilde{m}_\lambda + \sum_{\mu<\lambda}
c_{\lambda\mu}\tilde{m}_\mu,\quad
(ii)\;\langle P_\lambda^{HO},\tilde{m}_\mu\rangle_{HO} = 0\;
\text{for}\;\mu<\lambda.
\end{equation}
It is proved in \cite[p.\ 18]{hs:harmonic} that, for any $k\in V_{HO}$, the
$P^{HO}_\lambda(x;k)$ are mutually orthogonal with respect to the inner product
\eqref{e:HOorth}.

In the following theorem the zonal spherical functions 
on the symmetric space $U/K$ are expressed as $BC$ type Heckman-Opdam polynomials
(cf.\ \cite[p.\ 76]{hs:harmonic}).
\begin{thm}\label{th:cl-zonal}
Under restriction to the maximal torus $\TT$, the zonal $K'$-spherical function
$\varphi(\lambda)$ \textup{(}$\lambda\in P^+_K$\textup{)} is
mapped onto \textup{(}a scalar multiple of\textup{)} 
the Heckman-Opdam hyper\-geo\-met\-ric polynomial
$P^{HO}_{\lambda^\natural}(x;k)$ 
with $k_i=\frac{1}{2}m_i$ \textup{(}$i=1,2,3$\textup{)}.
\end{thm}
Of course, the zonal $K$-spherical functions can be described in the 
same way, since the subgroups $K$ and $K'$ are conjugate. 

For later purposes, it is convenient to rewrite the zonal spherical functions
in terms of generalized Jacobi polynomials, which are defined as follows.
Write $\CC[x]^{\goS}$ for the algebra of symmetric polynomials in
the variables $x_i$ ($1\leq i\leq l$). A linear basis of this algebra
is formed by the monomial symmetric polynomials 
$m_{\lambda}(x):=\sum_{\mu\in \goS\lambda}x^{\mu}$ ($\lambda\in P_\Sigma^+$).
Let $V_J$ denote the set of triples $(\al,\be,\tau)$ of real numbers
such that $\al,\be>-1$ and $\tau>0$. For any $(\al,\beta,\tau)\in V_J$,
we define an inner product 
$\langle \cdot,\cdot\rangle_J = \langle\cdot,\cdot\rangle_J^{\al,\be,\tau}$
on $\CC[x]^{\goS}$ by
\begin{equation}\label{e:Jorth}
\langle P,Q\rangle_J := \int_{x_1=0}^1 \cdots \int_{x_l=0}^1 
P(x) \overline{Q(x)} \Delta_J(x;\al,\be;\tau) dx,
\end{equation}
with $dx=dx_1\ldots dx_n$ and 
\begin{equation}
\Delta_J(x; \al,\be;\tau) := \prod_{i=1}^l (1-x_i)^\be x_i^\al |\Delta(x)|^{2\tau},
\end{equation}
where $\Delta(x):=\prod_{i<j} (x_i-x_j)$ is the Vandermonde determinant.
The generalized Jacobi polynomials 
$\lbrace P^J_\lambda(x;\al,\be;\tau) \mid \lambda\in P^+_\Sigma\rbrace$
(cf.\ Vretare \cite{vr:jacobi}) can now be defined by the same 
conditions \eqref{e:defpolHO} with $\tilde{m_\lambda}$ replaced by $m_\lambda$ and 
$\langle\cdot ,\cdot \rangle_{HO}$
replaced by $\langle\cdot ,\cdot \rangle_J$. 

As shown in \cite[\S3]{stok-koor:limit}, the simple change of variables 
$x_i\mapsto -\frac{1}{4}(x_i+x_i^{-1} -2)$ transforms 
the symmetric polynomials $P^J_\lambda$ into (a scalar multiple of)
the $\FSW$-invariant Laurent polynomials $P^{HO}_\lambda$.
The parameter correspondence is
\begin{equation*}
\al = k_1+k_3 -\frac{1}{2}, \quad \be= k_3-\frac{1}{2},
\quad \tau=k_2.
\end{equation*}
Under this change of variables, the orthogonality relations
of the Heckman-Opdam polynomials with respect to $\langle \cdot,\cdot\rangle_{HO}$ 
become the orthogonality relations of the generalized 
Jacobi polynomials $P^J_\lambda$ with respect to $\langle \cdot,\cdot\rangle_J$. 
The integral of the weight function $\Delta_J(x;\al,\be;\tau)$ was
first evaluated by Selberg in his well-known paper \cite{sb:integral}:
\begin{equation}\label{e:selberg}
\begin{split}
\int_{x_1=0}^1 \cdots \int_{x_l=0}^1 &\prod_{i<j} |x_i-x_j|^{2\tau}
\prod_{i=1}^l x_i^\al  (1-x_i)^\be dx = \\
&\prod_{j=1}^l {\Ga(\al+1+(j-1)\tau)\Ga(\be+1+(j-1)\tau)\Ga(j\tau+1)\over
\Ga(\al+\be+2+(l+j-2)\tau)\Ga(\tau+1)}.
\end{split}
\end{equation}
By the above-described relation between BC type Heckman-Opdam polynomials
and generalized Jacobi polynomials,   
the zonal spherical functions in Theorem \ref{th:cl-zonal}
can be rewritten as generalized Jacobi polynomials
with parameter values $\al = n-2l$, $\be=0$ and $\tau=1$.

Observe that the zonal spherical
functions, being matrix coefficients of irreducible representations, 
are mutually orthogonal with respect to the $L^2$ inner product
on $A\subset L^2(U,dg)$, where $dg$ is the normalized Haar measure on $U$ 
(Schur orthogonality). The restriction of the $L^2$ inner product to the algebra $\FSH$ 
of bi-$K'$-invariant matrix coefficients coincides under the
isomorphism of Theorem \ref{th:cl-Winvar} 
with the inner product $\langle .,. \rangle_{HO}$ on ${\mathbb{C}}[x^{\pm}]^{\FSW}$
up to a non-zero positive constant. This constant can
be explicitly determined using the evaluation of the Selberg integral \eqref{e:selberg}.
%%%%%%%%%%%%%%%%%%%%%%%%%%%%%%%%%%%%%%%%%%%%%%%%%%%%%%%%%%%%%%%%%%%%%%%%%%%%%
%%
%%  SECTION 3  BC TYPE Q-HYPERGEOMETRIC ORTHOGONAL POLYNOMIALS
%%
%%
%%%%%%%%%%%%%%%%%%%%%%%%%%%%%%%%%%%%%%%%%%%%%%%%%%%%%%%%%%%%%%%%%%%%%%%%%%%%%
\section{$BC$ type $q$-hypergeometric orthogonal polynomials}
\label{section:polynomials}
In this section we recall three families of multivariable ($BC$ type)
basic hypergeometric orthogonal polynomials. Later on, certain subfamilies 
of these multivariable orthogonal polynomials will be interpreted
as zonal spherical functions on quantizations of the classical 
complex Grassmannian. These families are $q$-analogues of the generalized
Jacobi polynomials (or, equivalently, $q$-analogues 
of the BC type Heckman-Opdam polynomials) and multivariable analogues
of well-known families of basic hypergeometric orthogonal polynomials
occurring in the Askey tableau \cite{aw:memoir}.

We start with introducing some notations and conventions.
We fix $0<q<1$ in the remainder of this paper. 
The q-shifted factorial is given by
\begin{equation*}
(a;q)_i := \prod_{j=0}^{i-1} ( 1-q^ja) \quad 
(i\in\ZZ_+),\quad
(a;q)_{\infty} := \prod_{j=0}^{\infty}(1-q^ja).
\end{equation*}
We use the shorthand notation
\begin{equation*}
(a_1,\ldots,a_m;q)_i:=\prod_{j=1}^m(a_j;q)_i 
\quad (i\in\ZZ_+\cup \{\infty\})
\end{equation*}
for products of $q$-shifted factorials.
Fix $N\in\ZZ$ and $\alpha,\beta\in\CC$.
The Jackson $q$-integral truncated at $N$ is defined by 
\begin{gather*}
\int_{\alpha}^{\beta}f(x)d_{q,N}x :=
 \int_{0}^{\beta}f(x)d_{q,N}x-\int_{0}^{\alpha}f(x)d_{q,N}x,\notag\\
\int_{0}^{\beta}f(x)d_{q,N}x := \sum_{k=0}^Nf(\beta q^k)(\beta q^k-\beta 
q^{k+1})\quad\text{if}\quad N\geq 0,\\
\int_{\alpha}^{\beta}f(x)d_{q,N}x := 0\quad\text{if}\quad N<0.\notag
\end{gather*}
The (non-truncated) Jackson $q$-integral is defined by
\begin{equation*}
\int_{\alpha}^{\beta}f(x)d_qx := \lim_{N\to \infty} 
\int_{\alpha}^{\beta}f(x)d_{q,N}x,
\end{equation*}
provided the limit exists.
The three families of multivariable orthogonal 
polynomials we recall in this section are the  
Koornwinder polynomials (cf.\ \cite{koor:maw}, \cite{stok:aw-orth},
\cite{stok:residue}),
and the multivariable big and little $q$-Jacobi polynomials 
(cf.\ \cite{stok:jac}, \cite{stok:residue}). 
The Koornwinder polynomials are multivariable
analogues of the well-known Askey-Wilson polynomials.
In the multivariable setting, these families depend on one additional 
deformation parameter $t\in (0,1)$. 
Since we only need the case $t=q^k$ with $k\geq 1$ integer
in this paper, we restrict ourselves to recalling the definitions and results
for this special case.
We first specify the parameter domain $V_X$ and the orthogonality inner product
$\langle \cdot, \cdot \rangle_X$ for each family ($X=K, B, L$), where
$K$ stands for Koornwinder polynomials, and $B$ resp.\ $L$ stands
for big resp.\ little $q$-Jacobi polynomials.

For the case $X=K$ we refer to \cite{koor:maw}, 
\cite{stok:aw-orth}, \cite{stok:residue}.
Take $V_{K}$ to be the set of quadruples  $\underline{t}:=(t_0,t_1,t_2,t_3)$ such that 
\begin{itemize}
\item[(1)] $t_0,t_1,t_2,t_3$ are real, or appear in complex conjugate pairs,
\item[(2)] $t_it_j \not\in \RR_{\geq 1}$ for $0\leq i<j\leq 3$.
\end{itemize}
Fix $\underline{t}\in V_{K}$. For $e\in \lbrace t_0,t_1,t_2,t_3\rbrace$ with $|e|>1$, 
let $N_e\in\ZZ$ be the 
largest integer such that $|eq^{N_e}|>1$. Take $N_e:=-1$ if $|e|\leq 1$.
Let $\int_{\TT_k}dx_1 \cdots dx_k$  denote $k$-fold contour integration along the
complex unit circle in counterclockwise direction.
For $0\leq m\leq l$ define a sesquilinear form 
$\langle \cdot, \cdot\rangle_{m,q,t}^{\underline{t}}$ 
on the space $\CC[x^{\pm 1}]^\FSW$ by
\begin{multline*}
\langle P, Q\rangle_{m} := \frac{2^m\binom{l}{m}}{(2\pi i)^{l-m}} 
\sum_{e_1,\ldots,e_m} \int_{x_1=0}^{e_1}\cdots
\int_{x_m=0}^{e_m}\int\cdots\int_{(x_{m+1},\ldots,x_l)\in \TT_{l-m}}\\
P(x)\overline{Q(x)}\Delta_{m}(x)
\frac{d_{q,N_{e_1}}x_1}{(1-q)x_1}\ldots \frac{d_{q,N_{e_m}}x_m}{(1-q)x_m}
\frac{dx_{m+1}}{x_{m+1}}\ldots \frac{dx_l}{x_l}.
\end{multline*}
Here the sum is taken over all $e_j\in\lbrace t_0,t_1,t_2,t_3\rbrace$, and
the partial weight function 
$\Delta_{m}(x;\underline{t};q,q^k)$ is defined by
\begin{equation*}\label{e:weight}
\Delta_{m}(x):= \prod_{r=1}^m w_{1}(x_r;e_r;f_r,g_r,h_r;q)
\prod_{s=m+1}^l w_{2}(x_s;\underline{t};q) 
\prod_{\stackrel{{\scriptstyle{i<j}}}{{\scriptstyle{\ep,\ep'\in \{\pm
1\}}}}}(x_i^{\ep}x_j^{\ep'};q)_k,
\end{equation*}
with $w_2$ the continuous weight function of the one-variable Askey-Wilson polynomials,
\begin{equation*}
w_2(x; \underline{t};q) := \frac{(x^2, x^{-2};q)_\infty}
{(t_0x, t_0/x, t_1x, t_1/x, t_2x, t_2/x, t_3x, t_3/x;q)_\infty},
\end{equation*}
and 
\begin{equation*}
w_1(eq^i;e;f,g,h;q):=\hbox{res}_{x=eq^i}\Bigl(\frac{w_1(x;\underline{t};q)}{x}\Bigr),
\end{equation*}
where the $f,g,h$ are such that $(e,f,g,h)$ is 
a permutation of $(t_0,t_1,t_2,t_3)$.
Then 
\begin{equation}\label{e:AWorth}
\langle P, Q\rangle_{K} := \sum_{m=0}^n\langle P,Q\rangle_m \quad 
(P,Q\in \CC[x^{\pm 1}]^\FSW)
\end{equation}
defines a positive definite 
inner product $\langle\cdot,\cdot\rangle_{K}=
\langle\cdot,\cdot\rangle_{K,q,t}^{\underline{t}}$, which is symmetric
in the parameters $t_0,t_1,t_2,t_3$.

The integral of the weight function $\FSI^{\underline{t}}_{K,q,t} := 
\langle 1, 1\rangle^{\underline{t}}_{K,q,t}$ is a $q$-analogue of Selberg's
integral \eqref{e:selberg}:
\begin{equation}\label{e:AWintegral}
\FSI^{\underline{t}}_{K,q,t} = 
2^l\,l!\,\prod_{j=1}^l \frac{(t,t^{l+j-2}t_0t_1t_2t_3;q)_\infty}
{(t^j,q,t_0t_1t^{j-1}, t_0t_2t^{j-1}, t_0t_3t^{j-1}, t_1t_2t^{j-1},
t_1t_3t^{j-1}, t_2t_3t^{j-1};q)_\infty}.
\end{equation}
For the case $|t_0|, |t_1|, |t_2|, |t_3| <1$ this was first proved by
Gustafson \cite{gus:selberg}. The general case was proved
in \cite{stok:residue}. 

For the multivariable big and little $q$-Jacobi polynomials 
($X=B,L$) we refer to \cite{stok:jac} and \cite{stok:residue}.
Take $V_B$ to be the set of quadruples $(a,b,c,d)$ such that $c,d>0$ and
\begin{equation}
a \in \bigl(-c/dq,1/q\bigr),\,\,  
b \in \bigl(-d/cq,1/q\bigr),
\end{equation}
or $a=cz, b=-d\bar{z} \text{ with } z \in \CC\setminus \RR$.
Fix $(a,b,c,d)\in V_B$ and define a positive definite inner product
$\langle\cdot,\cdot\rangle_B=
\langle \cdot,\cdot\rangle_{B,q,t}^{a,b,c,d}$ on $\CC[x]^{\goS}$ by
\begin{equation}\label{e:BJorth}
\langle P, Q  \rangle_{B} = \int_{x_1=-d}^c\int_{x_2=-d}^{c}\cdots
\int_{x_l=-d}^{c}P(x)\overline{Q(x)}\Delta_B(x)d_qx_1\ldots d_qx_l,
\end{equation}
with weight function $\Delta_B(x;a,b,c,d;q,q^k)$ given  by
\begin{equation}\label{e:DeltaBJ}
\Delta_{B}(x) = \Delta(x) \prod_{i=1}^l w_B(x_i)
\prod_{i<j}x_i^{2k-1}\bigl(q^{1-k}x_j/x_i;q\bigr)_{2k-1}.
\end{equation}
Here  $\Delta(x):=\prod_{1\leq i<j\leq l}(x_i-x_j)$ is 
the Vandermonde determinant and
\begin{equation*}
w_B(x;a,b,c,d;q):= \frac{(qx/c,-qx/d;q)_{\infty}}
{(qax/c,-qbx/d;q)_{\infty}}
\end{equation*}
is the weight function of the one-variable big $q$-Jacobi polynomials.
The evaluation of the integral $\FSI^{a,b,c,d}_{B,q,t}:= 
\langle 1, 1\rangle^{a,b,c,d}_{B,q,t}$ was first conjectured by Askey
\cite{ask:selberg} and subsequently proved by
Evans \cite{evans:beta-gamma}:
\begin{equation}\label{e:integralBJ}
\begin{split}
\FSI^{a,b,c,d}_{B,q,t} = 
l!q^{k^2{l\choose 3} -{k\choose 2}{l\choose 2}}&\prod_{i=1}^l
\frac{\Ga_q(\al +1 +(i-1)k)\Ga_q(\be + 1 +(i-1)k)\Ga_q(ik)}
{\Ga_q(\al + \be +2 +(l+i-2)k)\Ga_q(k)}\\
\times \prod_{i=1}^l &\frac{(-d/c,-c/d;q)_\infty (cd)^{1+(i-1)k}}
{(-q^{\al+1+(i-1)k}d/c, -q^{\be +1 +(i-1)k}c/d;q)_\infty
(c+d)},
\end{split}
\end{equation}
where $t=q^k$ ($k\geq 1$ integer), $a=q^\al$, $b=q^\be$ and with
\[
\Gamma_q(a):=(1-q)^{1-a}\frac{\bigl(q;q\bigr)_{\infty}}{\bigl(q^a;q\bigr)_{\infty}}
\quad (a\notin -\ZZ_+)
\]
the $q$-Gamma function. Note that \eqref{e:integralBJ} is another
$q$-analogue of Selberg's integral \eqref{e:selberg}.

Finally, take $V_L$ to be the set of pairs $(a,b)$ such that
$a\in (0,1/q)$ and $b\in (-\infty, 1/q)$. Fix
$(a,b)\in V_L$ and define a positive definite inner product 
$\langle\cdot,\cdot\rangle_L=
\langle\cdot,\cdot\rangle_{L,q,t}^{a,b}$ on $\CC[x]^{\goS}$ by
\begin{equation}\label{e:LJorth}
\langle P, Q \rangle_{L} = \int_{x_1=0}^1\ldots\int_{x_l=0}^1 P(x)
\overline{Q(x)}\Delta_L(x)d_qx_1\ldots d_qx_l,
\end{equation}
with weight function $\Delta_L(x;a,b;q,q^k)$ given by
\begin{equation}\label{e:DeltaLJ}
\Delta_L(x) = \Delta(x) \prod_{i=1}^l w_L(x_i)
\prod_{i<j} x_i^{2k-1}\bigl(q^{1-k}x_j/x_i;q\bigr)_{2k-1},
\end{equation}
and
\begin{equation*}
w_L(x;a,b;q):=\frac{(qx;q)_{\infty}}
{(qbx;q)_{\infty}}x^{\alpha}\quad (a=q^{\alpha})
\end{equation*}
the weight function of the one-variable little $q$-Jacobi polynomials.
The evaluation of the integral $\FSI^{a,b}_{L,q,t}:= 
\langle 1, 1\rangle^{a,b}_{L,q,t}$
was conjectured by Askey \cite{ask:selberg} and proved independently
by Habsieger \cite{hab:selberg} and Kadell \cite{kad:selberg}:
\begin{equation}\label{e:integralLJ}
\begin{split}
\FSI^{a,b}_{L,q,t}=  l!q^{k(\al +1){l\choose 2} + 2k^2 {l\choose 3}}&{ }\\
\times \prod_{i=1}^l
&\frac{\Ga_q(\al +1 +(i-1)k)\Ga_q(\be + 1 +(i-1)k)\Ga_q(ik)}
{\Ga_q(\al + \be +2 +(l+i-2)k)\Ga_q(k)}, 
\end{split}
\end{equation}
where $t=q^k$ ($k\geq 1$ integer), $a=q^\al$, $b=q^\be$.
This integral is yet another $q$-analogue of Selberg's integral \eqref{e:selberg}.

We now  define the corresponding families of orthogonal polynomials. 

It will be convenient to write $m_{\lambda}^X$ ($\lambda\in P_\Sigma^+$) 
for $\tilde{m}_{\lambda}$ ($X=K$) or $m_{\lambda}$ ($X=B, L$).
Fix $t=q^k$ ($k\geq 1$ integer) and fix parameters in $V_X$.
The Koornwinder, 
multivariable big resp.\ little $q$-Jacobi polynomials
$\lbrace P_{\lambda}^X \mid \lambda\in P_\Sigma^+ \rbrace$ ($X=K,B$ resp.\ $L$)
are uniquely determined by the following two conditions:
\begin{equation}\label{e:defpol}
(i)\; P_\lambda^X = m_\lambda^X + \sum_{\mu<\lambda}
c_{\lambda\mu}^Xm_\mu^X,\quad
(ii)\;\langle P_\lambda^X,m_\mu^X\rangle_X = 0\;
\text{for}\;\mu<\lambda.
\end{equation}

After suitable rescaling of the parameters,
the Koornwinder polynomials tend to the Heckman-Opdam
hypergeometric polynomials associated with $BC_l$ 
in the limit $q\to 1$ (cf.\ \cite[\S4]{dj:polynomial}),
and the multivariable big resp.\ little $q$-Jacobi polynomials
tend to the generalized Jacobi polynomials in the limit $q\to 1$
(cf.\ \cite[Thm.\ 5.1]{stok-koor:limit}). 

The $P_\lambda^X$ are joint eigenfunctions of 
a certain second-order q-difference operator $D_X$ which is self-adjoint with
respect to $\langle \cdot, \cdot \rangle_X$:
\begin{equation}\label{e:difference-op}
D_X :=
\sum_{j=1}^n \phi_{X,j}^+(x)(T_{q,j}-\text{
Id}) + \phi_{X,j}^-(x)(T_{q^{-1},j}-\text{ Id})
 \quad (X=K,B,L), 
\end{equation}
where $T_{q^{\pm 1},j}$ is the (multiplicative) $q^{\pm 1}$-shift 
in the variable $x_j$. In the case $X=K$ one has explicitly
\begin{gather*}
\phi_{K,j}^+(x):=
\frac{(1-t_0x_j)(1-t_1x_j)(1-t_2x_j)(1-t_3x_j)}{(1-x_j^2)(1-qx_j^2)}
\prod_{i\neq  j}
\frac{(1-tx_ix_j)(1-tx_i^{-1}x_j)}{(1-x_ix_j)(1-x_i^{-1}x_j)},\\
\phi_{K,j}^-(x):=
\frac{(t_0-x_j)(t_1-x_j)(t_2-x_j)(t_3-x_j)}{(1-x_j^2)(q-x_j^2)}
\prod_{i\neq  j}\frac{(t-x_ix_j)(t-x_i^{-1}x_j)}
{(1-x_ix_j)(1-x_i^{-1}x_j)}.
\end{gather*} 
The eigenvalues in this case are given by
\begin{equation*}\label{e:AWeigenvalue}
E^{K}_{\lambda}(\underline{t};q,t) = 
\sum_{j=1}^l\left (q^{-1}t_0t_1t_2t_3t^{2l-j-1}
(q^{\lambda_j}-1)+t^{j-1}(q^{-\lambda_j}-1)\right ).
\end{equation*}
For the $q$-difference operators $D_X$ and their eigenvalues in the
cases $X=B,L$ see \cite{stok:jac}.

As shown in \cite{stok:aw-orth} ($X=K$) and \cite{stok:jac} ($X=B,L$),
the self-adjointness of $D_X$ with respect to $\langle\cdot,\cdot\rangle_X$
may be used to prove full orthogonality:
\begin{equation}\label{e:fullorth}
\langle P_\lambda^X,P_\mu^X \rangle_X = 0\quad (\lambda,\mu\in P_\Sigma^+, 
\lambda\neq \mu).
\end{equation}
We remark that for $\underline{t}\in V_{K}$ with
$|t_0|, |t_1|, |t_2|, |t_3| <1$, the orthogonality measure for the Koornwinder
polynomials reduces to the completely continuous orthogonality measure
which was considered for the first time by Koornwinder \cite{koor:maw}.

For certain values of the parameters, for instance $t_0t_1t_2t_3\in [-q,1)$,
one has  $E_\lambda^{K}\neq  E_\mu^{K}$ if $\lambda<\mu$ (cf.\ 
\cite[Prop.\ 4.6]{stok-koor:limit}). The Koornwinder
polynomials may then be characterized by the conditions
\begin{equation}\label{e:AWqdiff}
\text{(i)}\;\; P^{K}_\lambda= \tilde{m}_\lambda + 
\sum_{\mu < \lambda}c_{\lambda\mu} \tilde{m}_\mu,\quad
\text{(ii)}\;\; D_{K}P^{K}_\lambda = E_\lambda^{K} P^{K}_\lambda.
\end{equation}
Using this characterization, we can read off 
certain elementary symmetry properties of the Koornwinder
polynomials from the corresponding symmetry properties of Koornwinder's
second-order $q$-difference operator and its eigenvalues. In the following lemma we 
formulate two of them.
\begin{lem}\label{symmetryproperties}
Let $\underline{t}\in V_K$ such that $t_0t_1t_2t_3\in [-q,1)$. 
Then the Koornwinder polynomial $P_{\lambda}(.;\underline{t};q,t)$ is symmetric
in the four parameters $t_0,t_1,t_2$ and $t_3$, and satisfies
\[P_{\lambda}(x;-\underline{t};q,t)=
(-1)^{|\lambda|}P_{\lambda}(-x;\underline{t};q,t),\]
where $-x:=(-x_1,\ldots,-x_l)$ and, similarly, $-\underline{t}=(-t_0,-t_1,-t_2,-t_3)$.
\end{lem}
Write $|\lambda|:=\sum_{i=1}^n\lambda_i$ for $\lambda\in P^+_\Sigma$, and
$cx:=(cx_1,\ldots,cx_l)$ for $c\in\CC$.
We have the following limit transitions from 
Koornwinder polynomials to multivariable big resp.\  
little q-Jacobi polynomials (cf. \cite{stok-koor:limit}, \cite{stok:limit}).
Fix $\lambda\in P_\Sigma^+$ and $k\geq 1$. For $(a,b,c,d)\in V_B$ we have the limit
transition
\begin{equation}\label{e:limitBJ}
\lim_{\ep \downarrow 0} 
\left (\frac{\ep (cd)^\frac{1}{2}}{q^\frac{1}{2}}\right)^{|\lambda |}
P^{K}_{\lambda}\Bigl(\frac{q^\frac{1}{2}x}
{\ep (cd)^\frac{1}{2}};\underline{t}_B(\ep);q,q^k\Bigr) = 
P_{\lambda}^B(x;a,b,c,d;q,q^k),
\end{equation}
where
\begin{equation}\label{tB}
\underline{t}_B(\ep):=\bigl(\ep^{-1}(qc/d)^\frac{1}{2},-\ep^{-1}
(qd/c)^\frac{1}{2},\ep a(qd/c)^\frac{1}{2},-\ep b(qc/d)^\frac{1}{2}\bigr).
\end{equation}
For $(a,b)\in V_L$ we have the limit transition
\begin{equation}\label{e:limitLJ}
\lim_{\ep \downarrow 0} 
\left(\frac{\ep}{q^\frac{1}{2}}\right)^{|\lambda |}
P^{K}_{\lambda}\Bigl(\frac{q^\frac{1}{2}x}{\ep};
\underline{t}_L(\ep);q,q^k\Bigr) =
P_{\lambda}^{L}(x;a,b;q,q^k),
\end{equation} 
where
\begin{equation}\label{tL}
\underline{t}_L(\epsilon):=\bigl(\ep^{-1}q^\frac{1}{2},-aq^\frac{1}{2},
\ep bq^\frac{1}{2},-q^\frac{1}{2}\bigr).
\end{equation}

To understand how $\FSW$-invariant Laurent polynomials can tend to
symmetric polynomials observe that 
$\lim_{\ep \downarrow 0}\ep^{|\lambda |}
\tilde{m}_\lambda(x/\ep) = m_\lambda(x)$ ($\lambda\in P_\Sigma^+$).

Let us write $m_i^X$ ($1\leq i\leq l$) for the orbit sum 
in $\CC[x^{\pm 1}]^\FSW$ 
corresponding to the $i$-th fundamental weight ($X=K$) resp.\ 
for the $i$-th elementary symmetric polynomial in $\CC[x]^\goS$ ($X=B,L$).
It is known (see \cite[Ch.~VI, \S4, Thm.\ 1]{bou:racines})
that the $m_i^X$ ($1\leq i\leq l$)  are algebraically independent and
generate the algebra $\CC[x^{\pm 1}]^\FSW$ resp.\ $\CC[x]^\goS$.
It therefore makes sense to introduce the notation 
\begin{equation}\label{altnotat}
\widehat{P}^X_\lambda(m_1^X(x), \ldots, m_l^X(x)) :=
P^X_\lambda(x)\quad (\lambda\in P^+_\Sigma).
\end{equation}
The $\widehat{P}^X_\lambda$ are (non-symmetric) polynomials in $l$ variables.
We may reformulate the above limit transitions in terms of these polynomials
in the following way.
\begin{equation}\label{e:limitBJ-rewritten}
\begin{split}
\lim_{\ep \downarrow 0} (s_\ep^{-1})^{|\lambda |}
\widehat{P}^{K}_{\lambda} (s_\ep y_1, \ldots, s_\ep^l y_l;
&\underline{t}_B(\ep);q,q^k) =\\ 
&\widehat{P}_{\lambda}^B(y_1, \ldots, y_l;a,b,c,d;q,q^k)
\end{split}
\end{equation}
with $s_\ep:= q^{\frac{1}{2}}/\ep(cd)^{\frac{1}{2}}$, and
\begin{equation}\label{e:limitLJ-rewritten}
\lim_{\ep \downarrow 0} (s_\ep^{-1})^{|\lambda |}
\widehat{P}^{K}_{\lambda}(s_\ep y_1, \ldots ,s_\ep^l y_l; 
\underline{t}_L(\ep);q,q^k) =
\widehat{P}_{\lambda}^{L}(y_1, \ldots , y_l;a,b;q,q^k)
\end{equation} 
with $s_\ep := q^\frac{1}{2}/\ep$.

Important for our applications of the multivariable orthogonal polynomial
theory to the study of zonal spherical functions on quantum Grassmannians
is the fact that the limit transitions \eqref{e:limitBJ} and \eqref{e:limitLJ}
extend to rigorous limits on the level of the orthogonality measure. 
This means that
the orthogonality measure of the Koornwinder polynomials
tends (in a suitable weak sense) to a non-zero multiple of 
the orthogonality measure of the multivariable 
big respectively little $q$-Jacobi polynomials
in the limit \eqref{e:limitBJ} resp.\ \eqref{e:limitLJ} (see \cite{stok:residue}).
In particular, the quadratic norms of the multivariable
big and little $q$-Jacobi polynomials can be derived by computing 
the limits \eqref{e:limitBJ} and \eqref{e:limitLJ} in the quadratic norm 
expressions of the Koornwinder polynomials (cf. \cite{stok:residue}).

We reformulate this observation as follows.
For parameters  $\underline{t}\in V_K$, $(a,b,c,d)\in V_B$ and $(a,b)\in V_L$,
let the renormalized quadratic norms
$N_K(\la):=N_K(\la;\underline{t};q,t)$, $N_B(\la):=N_B(\la;a,b,c,d;q,t)$
and $N_L(\la):=N_L(\la;a,b;q,t)$ for $\la\in P_{\Sigma}^+$ be defined by
\begin{equation}\label{NX}
N_X(\la):=\frac{\langle P_{\la}^X,P_{\la}^X\rangle_X}{\langle 1,1\rangle_X},\quad
(X=K,B,L).
\end{equation}
By the positive-definiteness of $\langle .,. \rangle_X$,
it follows that $N_X(\la)$ is strictly positive for all $\la\in P_{\Sigma}^+$. 
\begin{prop}\label{normlimit}\textup{(}\cite{stok:residue}\textup{)}
For $(a,b,c,d)\in V_B$, we have
\[
\lim_{\ep\downarrow 0}\,\bigl(\ep(cd/q)^{\frac{1}{2}}\bigr)^{2|\la|}
N_K(\la;t_B(\ep);q,t)=N_B(\la;a,b,c,d;q,t),\qquad (\la\in P_{\Sigma}^+).
\]
For $(a,b)\in V_L$, we have
\[
\lim_{\ep\downarrow 0}\,\bigl(\ep q^{-\frac{1}{2}}\bigr)^{2|\la|}N_K(\la;t_L(\ep);q,t)
=N_L(\la;a,b;q,t),\qquad (\la\in P_{\Sigma}^+).
\]
\end{prop}
\begin{rem}
For the explicit evaluations of the quadratic norms of 
the Koornwinder polynomials and the multivariable big and little
$q$-Jacobi polynomials in terms of products and quotients of $q$-Gamma functions,
we refer the reader to \cite{stok:residue}. The quadratic norm evaluations for
the Koornwinder polynomials were derived earlier in the paper \cite{vD1}. 
\end{rem}

%%%%%%%%%%%%%%%%%%%%%%%%%%%%%%%%%%%%%%%%%%%%%%%%%%%%%%%%%%%%%%%%%%%%%%%%%%%%%%%%%%%
%%
%%   SECTION 4    RECALL ON THE QUANTUM UNITARY GROUP
%%
%%
%%%%%%%%%%%%%%%%%%%%%%%%%%%%%%%%%%%%%%%%%%%%%%%%%%%%%%%%%%%%%%%%%%%%%%%%%%%%%%%%%%
\section{Preliminaries on the quantum unitary group}
\label{section:preliminaries}
Various aspects of the quantum unitary group have been studied in 
many different papers. Our main
references will be \cite{nym:gl} and \cite[\S1]{nou:macdonald}, which are based on
the $R$-matrix approach described in \cite{rtf:R}. 

The quantized coordinate ring $\Aq(\Mat(n,\CC))$ of the space of $n\times n$
complex matrices is defined as the  algebra with generators $t_{ij}$
($1\leq i,j\leq n$) and relations
\begin{gather*}
t_{ki}t_{kj}=qt_{kj}t_{ki},\quad t_{ik}t_{jk}=qt_{jk}t_{ik} \quad (i<j),\\
t_{il}t_{kj}=t_{kj}t_{il},
\quad t_{ij}t_{kl}-t_{kl}t_{ij}=(q-q^{-1})t_{il}t_{kj}
\quad (i<k,\,j<l).
\end{gather*}
In more compact notation, these relations may be written $RT_1T_2=T_2T_1R$.
Here $T_1:=T\ten I$, $T_2:=I\ten T$ ($I$ unit matrix, $\ten$ the Kronecker
product of matrices), and  
\begin{equation}\label{e:Rdef}
R := \sum_{ij} q^{\de_{ij}} e_{ii}\ten e_{jj} +
(q-q^{-1}) \sum_{i>j} e_{ij}\ten e_{ji},
\end{equation}
with the $e_{ij}$ ($1\leq i,j\leq n$) denoting the standard matrix units.
The matrix $R$ is invertible and satisfies the Quantum Yang-Baxter Equation.

The quantized coordinate ring $\AG$ of the general linear group
$G=GL(n,\CC)$ is defined by adjoining to $\Aq(\Mat(n,\CC))$ the
inverse $\detq^{-1}$ of the quantum determinant
\begin{equation*}
\detq:=\sum_{\sigma\in \goS_n}(-q)^{l(\sigma)}
t_{1\sigma(1)}\cdots t_{n\sigma(n)} \in \Aq(\Mat_q(n,\CC))
\end{equation*}
($l$ denoting the length function on $\goS$), which is central.
It follows from \cite[Lemma 1.5]{nym:gl} that $\Aq(\Mat(n;\CC))$ has no
zero divisors.

There is a unique Hopf algebra structure on $\AG$ such that $(t_{ij})$
becomes a matrix corepresentation.
The antipode $S\colon \AG\to\AG$ is given on the generators by
\begin{equation*}
S(t_{ij}):=(-q)^{i-j}\xi^{j^c}_{i^c}\detq^{-1}
\end{equation*}
with $i^c:=\lbrace 1,\ldots,n\rbrace\setminus \lbrace i \rbrace$, and 
with the quantum minor $\xi^I_J$ for subsets $I=\lbrace i_1<\ldots <i_r\rbrace, 
J=\lbrace j_1<\ldots <j_r\rbrace
\subset \lbrace 1,\ldots,n\rbrace$ defined by 
\begin{equation*}\label{e:qminor}
\xi^I_J:=\sum_{\sigma\in \goS_r}(-q)^{l(\sigma)}t_{i_1j_{\sigma(1)}}\cdots 
t_{i_rj_{\sigma(r)}}.
\end{equation*} 
$\AG$ becomes a Hopf $\ast$-algebra by requiring $(t_{ij})$ to be a unitary
matrix corepresentation. Im particular, this means that the $*$-structure
is given by $t_{ij}^\ast := S(t_{ji})$.
 We write $\Aq(U)=\Aq(U(n))$ for $\AG$ endowed
with this $\ast$-operation. The mapping $\tau:= \ast\circ S$ is a 
conjugate-linear involution on $\AU$ such that $\tau(t_{ij}) = t_{ji}$.

The quantized Borel subgroups 
$\Aq(B^\pm)$ of upper resp.\ lower triangular matrices
are defined as the Hopf quotients of $\AG$ by the relations
\begin{equation*}
t_{ij}=0\;(i>j)\quad \text{resp.}\quad t_{ij}=0 \;(i<j).
\end{equation*}
The corresponding projections will be denoted by
$\pi_\pm\colon \AG\to \Aq(B^\pm)$.
Note that the $z_i:=\pi_\pm(t_{ii})$ ($1\leq i\leq n$) 
in $\Aq(B^\pm)$ are invertible. Corresponding to the diagonal subgroup
$\TT\subset U(n)$ we have a natural surjective Hopf $\ast$-algebra morphism 
${}_{|\TT}$ of $\Aq(U)$ onto the Laurent polynomial algebra
$\AT:=\CC\lbrack z_1^{\pm 1},\ldots, z_n^{\pm 1} \rbrack$.

Next, we briefly recall the ``global'' description of finite-dimensional
corepresentations of $\AU$. For every $\lambda\in P$ 
(cf. \S \ref{section:classical}), we may define a
linear character (i.e.\ one-dimensional corepresentation) 
$z^{\lambda}:=z_1^{\lambda_1}\cdots z_n^{\lambda_n}$ of
$\Aq(B^\pm)$ or $\AT$. Using these linear characters it is completely straightforward
to define (highest) weight vectors in left or right $\AU$-comodules.
We take highest weight vectors of right resp.\ left 
$\AU$-comodules with respect to $\Aq(B^+)$ and $\Aq(B^-)$ respectively. For instance,
a highest weight vector of weight $\lambda\in P^+$ in a right
$\AU$-comodule $M$ is a non-zero vector $v\in M$ such
that 
\begin{equation*}
(\id\ten \pi_+)\circ \rho_M (v) = v\ten z^\lambda,
\end{equation*}
$\rho_M\colon M\to M\ten \AU$ denoting the comodule mapping.
Finite-dimensional $\AU$-comodules are then completely reducible and
unitarizable (see, for instance, \cite{nym:gl}, \cite{dijk-koor:cqg}). Recall
that a right $\AU$-comodule $M$ endowed with a positive definite inner
product (taken to be conjugate linear in the second variable) is
called unitary if
\begin{equation*}
\sum_{(v)(w)} \langle v_{(1)}, w_{(1)}\rangle w_{(2)}^\ast v_{(2)} =
\langle v,w\rangle 1 \quad (v,w\in M),
\end{equation*}
where the symbolic notation $\sum_{(v)} v_{(1)} \ten v_{(2)} := \rho_M(v)$
is used. The irreducible finite-dimensional $\AU$-comodules 
are parame\-trized by dominant weights
$\lambda\in P^+$ as in the classical case (cf. \cite[Th.\ 2.12]{nym:gl}). 
The irreducible right $\AU$-comodule with highest weight $\lambda\in P^+$
is denoted by $V_R(\lambda)$. The vector space $V_L(\lambda):= \Hom(V_R(\lambda),\CC)$
has a natural left $\AU$-comodule structure, which is also irreducible
of highest weight $\lambda$. If no confusion is
possible, we will write $V(\lambda)$ for the left comodule 
$V_L(\lambda)$ respectively for the right comodule $V_R(\lambda)$.
\begin{rem}\label{th:Kleftright}
Let $M$ be a finite-dimensional right $\AU$-comodule with comodule mapping
$\rho_M\colon M \to M \ten \AU$. Write $M^\circ$ for the vector
space complex conjugate to $M$ and $\si\colon M\ten \AU \to \AU\ten M$ for
the flip. Then the mapping 
\begin{equation}\label{e:conjugate}
\rho_M^\circ\colon M^\circ \to \AU \ten M^\circ, \quad
\rho_m^\circ := (\tau \ten \id)\circ \si\circ \rho_M,
\end{equation}
where $\tau=\ast\circ S$,
defines a left $\AU$-comodule structure on $M^\circ$.
In \eqref{e:conjugate} $\rho_M$ is considered as a conjugate 
linear map from $M^\circ$ to $M\otimes A_q(U)$,
and $\tau\ten \hbox{id}$ as a conjugate linear map from $A_q(U)\ten M$
to $A_q(U)\ten M^\circ$.
 
The assignment $M\mapsto M^\circ$ is a 1-1 correspondence between right and 
left $\AU$-comodules preserving weights and highest weights. Hence
$M^\circ$ is isomorphic to the left $\AU$-comodule $\Hom(M, \CC)$. 
A right $\AU$-comodule intertwiner $\Psi: M\rightarrow
N$ also intertwines the left $\AU$-comodule structures of $M^\circ$ and
$N^\circ$ (i.e.\ when $\Psi$ is considered as a map from $M^\circ$ to $N^\circ$).
\end{rem}

Recall that the comultiplication $\Delta\colon \AU \to \AU\ten\AU$ defines
a bicomodule structure on $\AU$.
Let $W(\lambda)\subset \AU$ ($\lambda\in P^+$) denote the subspace spanned
by the matrix coefficients of either $V_R(\lambda)$ or $V_L(\lambda)$. 
The irreducible decomposition of the bicomodule $\AU$ reads
\begin{equation}\label{e:peterweyl}
\AU = \bigoplus_{\lambda\in P^+} W(\lambda), \quad 
W(\lambda) \simeq V_L(\lambda) \ten V_R(\lambda).
\end{equation}
Let $h$ be the normalized Haar functional on $\AU$. It can be characterized
as the unique linear functional on $\AU$ which is zero on $W(\lambda)$
for $0\not=\lambda\in P^+$ and which sends $1\in \AU$ to $1\in {\mathbb{C}}$.
The subspaces $W(\lambda)$ are mutually orthogonal with respect to
the inner product $\langle \varphi, \psi\rangle := h(\psi^\ast\varphi)$.

We consider now in some more detail the vector corepresentation, its
dual representation, and their exterior powers.
Let $V$ denote the vector space $\CC^n$ with canonical basis $(v_i)$. 
$V$ becomes a right $\AU$-comodule (called vector corepresentation) with
\begin{equation}\label{e:Rvector}
\rho_V\colon V\mapsto V\otimes \AU,\quad
\rho_V(v_j):=\sum_{i=1}^nv_i\otimes t_{ij} \quad (1\leq j\leq n).
\end{equation}
$V$ is irreducible with highest weight $\tilde{\ep}_1$ and highest weight 
vector $v_1$. Note that the vectors $v_i$ have weight $\tilde{\ep}_i$.
The corepresentation $V$ is unitary with respect to the 
inner product $\langle v_i, v_j\rangle =\delta_{ij}$.

Let $V^\ast$ denote the linear dual of $V$ with dual basis $(v_i^\ast)$.
$V^\ast$ becomes a right $\AU$-comodule (contragredient of $V$) with
\begin{equation}\label{e:Rcovector}
\rho_{V^\ast}\colon V^\ast\mapsto V^\ast\otimes \AU,\quad
\rho_{V^\ast}(v_j^\ast):=\sum_{i=1}^nv_i^\ast\otimes t_{ij}^\ast\quad (1\leq j\leq n).
\end{equation}
$V^\ast$ is irreducible with highest weight $-\tilde{\ep}_n$ and highest weight 
vector $v_n^\ast$. Note that the vectors $v_i^\ast$ have weight $-\tilde{\ep}_i$. 
The corepresentation $V^\ast$ is unitary with respect to the inner
product $\langle v_i^\ast, v_j^\ast\rangle := q^{-\langle 2\rho, \tilde{\ep}_i\rangle}
\delta_{ij}$, where $\rho:= \sum_{k=1}^n (n-k)\tilde{\ep}_k$. This follows
{}from the well-known fact that $S^2(t_{ij}) = 
q^{\langle 2\rho, \tilde{\ep}_j-\tilde{\ep}_i\rangle}t_{ij}$ ($1\leq i,j\leq n$).

Let $\Lambda_q(V)$ resp.\ $\Lambda_q(V^\ast)$ denote the associative 
algebra generated by $v_1,\ldots,v_n$ resp.\ $v_1^\ast, \ldots, v_n^\ast$
with relations
\begin{equation}\label{e:extrel}
v_i\wedge v_i=0 \;(1\leq i\leq n),\quad v_i\wedge 
v_j=-q^{-1}v_j\wedge v_i \;(i<j)
\end{equation}
respectively 
\begin{equation}\label{e:extreldual}
v_i^\ast\wedge v_i^\ast=0\; (1\leq i\leq n),
\quad v_j^\ast\wedge v_i^\ast=-q^{-1}v_i^\ast\wedge
v_j^\ast \; (i<j).
\end{equation}
Then $\Lambda_q(V)$ resp.\ $\Lambda_q(V^\ast)$ 
inherits a natural right $\AU$-comodule structure from $V$ resp.\ $V^\ast$ 
by extending the comodule mapping $\rho_V$ resp.\ $\rho_{V^\ast}$
as a unital algebra homomorphism. 
$\Lambda_q(V)$ resp.\ $\Lambda_q(V^\ast)$ also has a natural grading
such that the generators $v_i$ resp.\ $v_i^\ast$ have degree $1$:
\begin{equation*}
\Lambda_q(V)=\bigoplus_{r=0}^n\Lambda_q^r(V),\quad 
\Lambda_q(V^\ast)=\bigoplus_{r=0}^n\Lambda_q^r(V^\ast).
\end{equation*}
Write $v_I:=v_{i_1}\wedge\cdots \wedge v_{i_r}$ resp.\ 
$v_I^\ast:=v_{i_r}^\ast\wedge\cdots \wedge v_{i_1}^\ast$
if $I=\lbrace i_1<\cdots <i_r\rbrace\subset \lbrace 1, \ldots, n\rbrace$.
Write $|I|$ for the cardinality of $I$.
Then the $v_I$ resp.\ $v_I^\ast$ ($|I|=r$) form 
a basis of $\Lambda_q^r(V)$ resp.\ $\Lambda_q^r(V^\ast)$.
One has the  multiplicative property
($I,J\subset \lbrace 1, \ldots, n\rbrace$)
\begin{equation*}\label{e:multiwedge}
v_I\wedge v_J=\sgn(I;J)v_{I\cup J},\quad 
v_I^\ast\wedge v_J^\ast:=\sgn\bigl(J;I\bigr)v_{I\cup J}^\ast
\end{equation*}
where 
\begin{equation*}
\sgn(I;J):=
\begin{cases} 0& \text{if $I\cap J\neq \emptyset$},\\ 
(-q)^{l(I;J)}& \text{if $I\cap J=\emptyset$},
\end{cases}
\end{equation*}
and $l(I;J):=|\lbrace (i,j)\in I\times J \, | \, i>j \rbrace|$.
The comodules $\Lambda_q^r(V)$ and $\Lambda_q^r(V^\ast)$ are irreducible
subcomodules of $\Lambda_q(V)$, and the coactions satisfy
\begin{equation}\label{e:action}
\rho_V(v_J)=\sum_{|I|=r}v_I\otimes \xi^I_J, \quad
\rho_{V^\ast}(v_J^\ast)=\sum_{|I|=r}v_I^\ast\otimes (\xi^I_J)^\ast \quad (|J|=r).
\end{equation}
{}From this it follows immediately that
\begin{equation}\label{e:deltaqminor}
\Delta(\xi^I_J)=\sum_{|K|=r}\xi^I_K\otimes\xi^K_J\quad (|I|,|J|=r)
\end{equation}
and hence
\begin{equation}\label{e:sumqminor}
\sum_{|K|=r}\xi^I_K S(\xi^K_J) = \delta_{I,J} \quad (|I|,|J|=r).
\end{equation}
We furthermore recall (cf.\ \cite[(3.2)]{nym:gl})
that the quantum minors satisfy
\begin{equation}\label{e:qminoranti}
(\xi^I_J)^\ast=S(\xi^J_I)=\frac{\sgn(J;J^c)}
{\sgn(I;I^c)}\xi^{I^c}_{J^c}\detq^{-1} \quad (|I|=|J|=r),
\end{equation} 
where $I^c:=\lbrace 1,\ldots,n\rbrace \setminus I$.

The $\AU$-comodule $\Lambda_q^r(V)$ resp.\ 
$\Lambda_q^r(V^\ast)$ has highest weight
$\Lambda_r:=\tilde\ep_1 + \cdots +\tilde\ep_r$ resp.\ 
$\Lambda_{n-r}-\Lambda_n = -\tilde\ep_{n-r+1} -\cdots -\tilde\ep_n$ with
highest weight vector $v_1\wedge \cdots \wedge v_r$ resp.\ 
$v_n^\ast \wedge \cdots \wedge v_{n-r+1}^\ast$.
The inner product on $\Lambda_q^r(V)$ such that
$\langle v_I, v_J\rangle = \delta_{I,J}$ is $\AU$-invariant.
On the space $\Lambda_q^r(V^\ast)$ we have the invariant
inner product $\langle v_I^\ast, v_J^\ast\rangle :=
\delta_{I,J} q^{-\langle 2\rho, \tilde{\ep}_I\rangle}$,
where $\tilde{\ep}_I:=\sum_{i\in I}\tilde{\ep}_{i}$.

Let $\Uq=U_q(\gog\gol(n,\CC))$ denote the quantized universal enveloping algebra
(cf.\ Drinfel'd \cite{drin:proc}, Jimbo \cite{jim:uqgl}) associated with
the Lie algebra $\gog=\gog\gol(n,\CC)$. In our 
notation we essentially adhere to Noumi \cite[\S 1]{nou:macdonald}. 
We refer to this last paper for any details not treated here. The
algebra $\Uq$ is generated by elements $q^h$ ($h\in P$) and $e_i, f_i$
($1\leq i\leq n-1$) subject to the quantized Weyl-Serre
relations. 

More useful for the purposes of this paper are
the $L$-operators $L^+_{ij},\,L^-_{ij}\in \Uq$ ($1\leq i,j \leq n$).
They generate $\Uq$ subject to certain commutation relations 
that may be conveniently expressed by means of the matrix $R$ defined in
\eqref{e:Rdef} (cf.\ \cite{rtf:R}). The matrices $L^\pm := (L^\pm_{ij})$
are upper resp.\ lower triangular, and $L^\pm_{ii} = q^{\pm \tilde{\ep}_i}$
($1\leq i\leq n$).
The Hopf $\ast$-algebra structure on $\Uq$ is uniquely determined by 
\begin{equation}\label{e:uqhopf}
\Delta(L^\pm_{ij}) = \sum_k L^\pm_{ik} \ten 
L^\pm_{kj}, \quad \ep(L^\pm_{ij}) = \delta_{ij},
\quad (L^\pm_{ij})^\ast = S(L^\mp_{ji}) \quad (1\leq i,j\leq n)
\end{equation}
The  involution
$\tau = \ast \circ S\colon \Uq \to \Uq$ acts on the generators as
\begin{equation}\label{e:Ltau}
\tau(L^\pm_{ij}) = L^\mp_{ji}\quad (1\leq i,j \leq n).
\end{equation}

There is a natural Hopf $\ast$-algebra 
duality $\langle \cdot\, , \, \cdot \rangle$
between $\Uq$ and $\AU$. This means in particular that we have 
\begin{equation}
\langle u, \varphi^\ast\rangle = \overline{\langle \tau(u), \varphi\rangle},
\quad
\langle u^\ast, \varphi\rangle = \overline{\langle u, \tau(\varphi)\rangle}
\quad (u\in \Uq, \varphi\in \AU).
\end{equation}
We write $\Uh$ for the subalgebra generated
by the $q^h$ ($h\in P$). It is Laurent polynomial
in the elements $q^{\tilde{\ep}_i}$ ($1\leq i\leq n$).
There is an induced Hopf $\ast$-algebra duality between
$\Uh$ and $\AT$ such that 
\begin{equation*}
\langle q^h\,,\, z^\lambda\rangle 
:= q^{\langle h, \lambda\rangle}, \quad
z^\lambda = z_1^{\lambda_1} \cdots z_n^{\lambda_n}\quad
(h, \lambda \in P).
\end{equation*}

For a right $\AU$-comodule $(M, \rho_M)$, 
the $\AU$-coaction $\rho_M$ can be ``differentiated'' 
using the Hopf algebra pairing $\langle\cdot, \cdot \rangle$.
This yields a left $\Uq$-module structure on $M$  (cf.\ \cite{nym:gl}).
To be precise, the left $\Uq$-action on $M$ is defined by
\begin{equation}\label{diffactionX}
X\cdot v:=\sum_{(v)}\langle X,v_{(2)}\rangle v_{(1)},\quad 
(X\in \Uq, v\in M),
\end{equation}
where $\rho_M(v)=:\sum_{(v)}v_{(1)}\ten v_{(2)}\in M\ten \AU$ for $v\in M$.
For example, differentiating the vector corepresentation
\eqref{e:Rvector}, the corresponding left $\Uq$-action gives rise to
an algebra homomorphism $\rho_V\colon
\Uq \to \End(V)$, which is uniquely determined by the formulas
\begin{equation}\label{e:Lvector}
R^\pm = \sum_{ij} e_{ij} \ten \rho_V(L^\pm_{ij}), \quad
(R^\pm)^{-1} = \sum_{ij} e_{ij} \ten \rho_V(S(L^\pm_{ij})).
\end{equation}
Here $R^-:=R^{-1}$, $R^+=PRP$, and $P=\sum_{i,j}e_{ij}\ten e_{ji}$ is
the permutation operator.
By differentiation of right $\AU$-coactions,
a 1-1 correspondence is obtained between
right $\AU$-comodule structures on a finite-dimensional vector space
$M$ and $P$-weighted left $\Uq$-module structures on $M$.
Recall that $M$ is $P$-weighted if it is spanned by vectors that transform under
$\Uh$ according to $q^h\cdot v = q^{\langle h, \lambda\rangle} v$
($\lambda\in P$). A highest weight vector $v$ of highest weight
$\lambda$ in a left
$\Uq$-module $M$ is then characterized by the conditions
$L^-_{ij}\cdot v = 0$ ($i>j$) (or, equivalently, $X_i^+\cdot v=0$ for $i\in [1,n-1]$) 
and $q^h\cdot v = q^{\langle h, \lambda\rangle} v$. 
There is a similar relationship between left $\AU$-comodules and right $\Uq$-modules. 
For a right $\Uq$-module $M$, a weight vector $0\not=v\in M$ is a highest weight
vector if $v\cdot L_{ij}^+=0$ ($i<j$) (or, equivalently, $v\cdot X_i^-=0$ for
$i\in [1,n-1]$).
\begin{rem}\label{starleftright}
For a right $\AU$-comodule $(M,\rho_M)$ and for an element $v\in M$, 
we write $v^\circ$ when considering $v$ as an  element in 
the left $\AU$-comodule $(M^\circ, \rho_M^\circ)$ (cf.\  Remark \ref{th:Kleftright}).
Then the differentiated right $\Uq$-module structure on $M^\circ$ is related
to the differentiated left $\Uq$-module structure on $M$ by
$v^\circ\cdot X=(X^*\cdot v)^\circ$, where $v\in M$ and  $X\in \Uq$.
\end{rem}
The coalgebra structure of $\AU$ naturally induces a 
$\AU$-bicomodule structure on $\AU$.
By differentiating this $\AU$-bicomodule structure, $\AU$ becomes a $\Uq$-bimodule 
with $\Uq$-symmetry. The action of the $L$-operators is then given by
\begin{equation}\label{e:Laction}
L_1^\pm \cdot T_2 = T_2 R^\pm, \quad T_2 \cdot L^\pm_1 = R^\pm T_2.
\end{equation}
Obviously, the irreducible decomposition of the $\Uq$-bimodule
$\AU$ is given by \eqref{e:peterweyl}. This decomposition may also
be characterized as the simultaneous eigenspace decomposition of
$\AU$ under the action of the center $\FSZ\subset \Uq$.

It can be shown that the pairing 
$\langle \cdot\, , \, \cdot \rangle$ is 
doubly non-degenerate (cf.\ \cite[Cor.\ 23, 54]{kli-schm:book}). 
In particular, $\AU$ can be embedded
as Hopf $*$-algebra into the Hopf $*$-algebra dual of $\Uq$.
The image under this embedding is the Hopf subalgebra spanned by matrix
elements of finite-dimensional $P$-weighted $\Uq$-modules.

%%%%%%%%%%%%%%%%%%%%%%%%%%%%%%%%%%%%%%%%%%%%%%%%%%%%%%%%%%%%%%%%%%%%%%%%%%%%%%%%%%%
%%
%%  SECTION 5   SPHERICAL COREPRESENTATIONS
%%
%%
%%%%%%%%%%%%%%%%%%%%%%%%%%%%%%%%%%%%%%%%%%%%%%%%%%%%%%%%%%%%%%%%%%%%%%%%%%%%%%%%%%
\section{Spherical corepresentations}\label{section:spherical}
We call $\AK:= \Aq(U(n-l))\ten \Aq(U(l))$ the quantized coordinate
ring of $K$. Corresponding to the embedding \eqref{e:Kembedding}
there is an obvious surjective Hopf $\ast$-algebra morphism
$\pi_K\colon \Aq(U(n)) \to \Aq(K)$. 
Write $A_q(U/K)$ for the right $A_q(K)$-fixed elements in $A_q(U)$, i.e.\
\begin{equation}\label{AUKeerst}
A_q(U/K):=\{ \varphi\in A_q(U) \, | \, (\hbox{id}\ten \pi_K)\circ \Delta(\varphi)=
\varphi\ten 1\}.
\end{equation}
Observe that $A_q(U/K)$ is a left $A_q(U)$-comodule $*$-subalgebra of $A_q(U)$.
The algebra $A_q(U/K)$ can be interpreted as a quantized algebra of functions
on the complex Grassmannian $U/K$.
 
For the study of $A_q(U/K)$ it is important to obtain explicit  
information about $\AK$-spherical corepresentations of
$\AU$, i.e.\  finite-dimensional right $\AU$-comodules with non-zero
$\AK$-fixed vectors.
 Recall that a vector $v$ in a right $\AU$-comodule
$M$ with comodule mapping $\rho_M\colon M\to M\ten \Aq(U)$ is $\AK$-fixed if
\begin{equation}
(\id\ten\pi_K)\circ \rho_M(v) = v \ten 1.
\end{equation}
One defines $\AK$-fixed vectors in left $\AU$-comodules in a similar way.
In this section we discuss the proof of the following theorem,
which states in particular that the pair $(\AU, \AK)$ is a quantum Gelfand pair:
\begin{thm}\label{th:gelfand}
Every finite-dimensional irreducible corepresentation of $\AU$ has at most one 
$\AK$-fixed vector \textup{(}up to scalar multiples\textup{)}. 
The finite-dimensional corepresentations with non-zero $\AK$-fixed vectors are
parametrized by the classical sublattice $P^+_K$ of spherical dominant weights
\textup{(}cf.\ \S \ref{section:classical}\textup{)}.
\end{thm}
\noindent
\begin{rem}\label{th:Kleftright2}
Let $M$ be a finite-dimensional right $\AU$-comodule.
It follows from Remark \ref{th:Kleftright} that
a vector $v\in M$ is $\AK$-fixed if and only if $v\in M^\circ$ 
is $\AK$-fixed. 
Hence, any statement about $\AK$-fixed vectors in right $\AU$-comodules 
immediately translates to a corresponding statement 
for left $\AU$-comodules and vice-versa.
\end{rem}
For the proof of Theorem \ref{th:gelfand} it suffices to show
that the irreducible decomposition of $V_R(\lambda)$ as
a right $A_q(K)$-comodule is the same as the decomposition of
the irreducible finite-dimensional representation of $U(n)$
with highest weight $\lambda$ when restricted to 
the subgroup $K$.
One way of establishing this result is by differentiating the coaction of $A_q(U)$ on
$V_R(\lambda)$ using the doubly non-degenerate Hopf algebra pairing between 
$A_q(U)$ and $U_q({\mathfrak{g}})$.
Then the desired result follows from well-known results on the representation
theory of quantized universal enveloping algebras.
This approach is quite general, and is treated in more detail
in \cite{SD}.

In this section we discuss another proof of Theorem \ref{th:gelfand} which does not
rely on the quantized universal enveloping algebra technique. The strategy will
be to relate the decomposition of the restriction to $\AK$ of
the right $\AU$-comodule $V_R(\lambda)$ ($\lambda\in P^+$) to characters
on the maximal torus $\TT$.
The following general result about 
corepresentation theory of semisimple coalgebras is needed (a coalgebra is
said to be semisimple if every finite-dimensional $A$-comodule is completely reducible).
\begin{prop}\label{semisimpel}
Let $A$ and $B$ be semisimple coalgebras.
Then every finite-dimen\-si\-onal $A\otimes B$-comodule
is completely reducible.
Write $\lbrace V_{\alpha} \, | \, \alpha\in \hat{A} \rbrace$ and
$\lbrace V_{\beta} \, | \, \beta\in \hat{B} \rbrace$
for a complete set of mutually inequivalent, irreducible, finite-dimensional
right $A$- and ${B}$-comodules, respectively. Then
\begin{equation}\label{compset}
\lbrace V_{\alpha}\boxtimes 
V_{\beta} \, | \, \alpha\in\hat{A}, \beta\in \hat{B} \rbrace
\end{equation}
is a complete set of mutually inequivalent irreducible finite-dimensional
right ${A}\otimes B$-comodules. Here $V_{\alpha}\boxtimes V_{\beta}=V_{\alpha}\ten
V_{\beta}$ as a vector space and it has right comodule structure given by
$\rho_{\alpha}\boxtimes \rho_{\beta}:=
\sigma_{23}\circ (\rho_{\alpha}\otimes \rho_{\beta})$, where
$\sigma_{23}$ is the flip of the second and third tensor component and where 
$\rho_{\alpha}$ and $\rho_{\beta}$ are the right comodule mappings of 
$V_{\alpha}$ and $V_{\beta}$, respectively.
\end{prop}
The proof of the proposition is similar to the analogous, well-known result 
for tensor products of semisimple algebras and will therefore be omitted here.

For $\lambda\in P^+$ with $\lambda_n\geq 0$, 
define the Schur polynomial $s_{\lambda}(z)\in\AT$ by
\begin{equation*}
s_\lambda(z):=\Delta^{-1}(z)\sum_{w\in \goS_n}(-1)^{l(w)}z^{w(\lambda+\rho)},
\end{equation*}
with $\Delta(z):=\prod_{i<j} (z_i-z_j)$ the Vandermonde determinant.
For arbitrary $\lambda\in P^+$ with $\lambda_n\geq -m$ \, ($m\in \ZZ$) 
define $s_\lambda(z):=z^{-m\Lambda_n}s_{\lambda+m\Lambda_n}(z)\in\AT$.
Then the $s_\lambda$ ($\lambda\in P^+$) are well-defined and form
a basis of the subalgebra $\AT^{\goS_n}$ of symmetric Laurent polynomials. 
Recall that the character of a finite-dimensional
$A_q(U)$-comodule $M$ is defined by 
$\chi_M := \sum_i \pi_{ii}\in \AU$, where the $\pi_{ij}\in A_q(U)$ are the matrix 
coefficients of $M$ with respect to a basis of $M$.
The character $\chi_M$ is independent  
of the particular choice of basis for $M$. 
As shown in \cite[(3.22)]{nym:gl},
the character $\chi_\lambda\in\AU$ of the irreducible comodule
$V_R(\lambda)$ satisfies
\begin{equation}\label{characterformula}
(\chi_{\lambda})_{|\TT}= s_{\lambda}(z) \quad (\lambda\in P^+),
\end{equation}
as in the classical case ($q=1$).
\begin{prop}\label{th:decomp} 
Let $\lambda\in P^+$. The restriction of the $\AU$-comodule
$V_R(\lambda)$ to $\AK$ decomposes as
\begin{equation}\label{e:decomp1}
V_R(\lambda)\simeq\bigoplus_{\mu,\nu}\,
(V_R(\mu)\boxtimes V_R(\nu))^{\oplus c_{\mu,\nu}^\lambda},
\end{equation}
the sum ranging over $\mu\in P_{n-l}^+$, $\nu\in P_l^+$.
Here the $c^\lambda_{\mu,\nu}$ are the non-negative integers
characterized by
\begin{equation}\label{e:decomp2}
s_\lambda(z_1,\ldots,z_n)=\sum_{\mu,\nu}
c_{\mu,\nu}^\lambda s_{\mu}(z_1,\ldots,z_{n-l})s_{\nu}(z_{n-l+1},\ldots,z_n),
\end{equation}
the sum ranging over $\mu\in P_{n-l}^+$, $\nu\in P_l^+$.
\end{prop}
\begin{proof}
There exists a decomposition \eqref{e:decomp1} for certain uniquely determined
non-negative integers $c_{\mu,\nu}^{\lambda}$ by the previous proposition.
It follows from \eqref{characterformula} that 
the $c_{\mu,\nu}^{\lambda}$ satisfy \eqref{e:decomp2}, since 
$\chi_{M\boxtimes N} = \chi_M\ten\chi_N\in \AK$ for a 
finite-dimensional right $A_q(U(n-l))$-comodule $M$ 
and a finite-dimensional right $A_q(U(l))$-comodule $N$.
\end{proof}
We conclude from Proposition \ref{th:decomp} 
that the abstract decomposition of an arbitrary finite-dimensional
right $\AU$-comodule $M$ into irreducible $\AK$-comodules is 
the same as in the classical ($q=1$) case. Hence, at this point we
see that Theorem \ref{th:gelfand} is a consequence of Theorem
\ref{th:cl-gelfand}.

\begin{rem} 
The proof of Theorem \ref{th:gelfand} can also be derived from
Proposition \ref{th:decomp} using the Littlewood-Richardson rule
(cf.\ Macdonald \cite{mac:symmetric}),
which is a combinatorial rule for computing the coefficients
$c^{\lambda}_{\mu,\nu}$ in \eqref{e:decomp2}.
\end{rem}

%%%%%%%%%%%%%%%%%%%%%%%%%%%%%%%%%%%%%%%%%%%%%%%%%%%%%%%%%%%%%%%%%%%%%%%%%%%%%%%%%%
%%
%%
%%   SECTION 6  A ONE-PARAMETER FAMILY OF QUANTUM GRASSMANNIANS
%%
%%
%%
%%%%%%%%%%%%%%%%%%%%%%%%%%%%%%%%%%%%%%%%%%%%%%%%%%%%%%%%%%%%%%%%%%%%%%%%%%%%%%%%%%
\section{A one-parameter family of quantum Grassmannians}
\label{section:nds}
In this section we define a family of quantum Grassmannians depending 
on one real parameter $-\infty < \si<\infty$ (cf.\ \cite[\S 2]{nds:gras}).
The key ingredient in the definition will be 
the $n\times n$ complex matrix $J^\sigma$ defined by
\begin{equation}\label{Jsi}
J^\sigma := \sum_{1\leq k\leq l} 
(1 - q^{2\sigma}) e_{kk} + \sum_{l<k<l'} e_{kk}
-\sum_{1\leq k\leq l} q^\sigma e_{kk'}-\sum_{1\leq k\leq l} q^\sigma e_{k'k},
\end{equation}
where $k':= n-k+1$ ($1\leq k\leq n$). Observe that
$\lim_{\si\to\infty} J^\si = J^\infty$, where $J^\infty$ is defined by
\begin{equation}\label{Jinfty}
J^{\infty}:=\sum_{k=1}^{n-l}e_{kk}.
\end{equation}
The subspace $\ksi\subset \Uq$
is by definition spanned by the coefficients of the matrix
\begin{equation}\label{e:ksidef}
L^+ J^\sigma - J^\sigma L^- \in \End(V) \ten \Uq.
\end{equation}
It follows from \eqref{e:uqhopf} 
that $\ksi$ is a two-sided coideal in $\Uq$, i.e.\ 
$\Delta(\ksi)\subset \Uq\ten \ksi + \ksi \ten \Uq$ and
$\ep(\ksi)=0$. This remains true when 
$J^{\si}$ is replaced by any $n\times n$ matrix $J$ in 
the definition of $\ksi$. Moreover, since $J^{\si}$ is a symmetric matrix, 
it follows from \eqref{e:Ltau} that $\ksi$ is $\tau$-invariant. 

Define the subalgebra $\Aksi\subset \AU$ as the subspace of
all left $\ksi$-invariant elements in $\AU$, i.e.\ all
$a\in \AU$ such that $\ksi\cdot a=0$. As is well-known (cf.\ for instance
\cite[Prop. 1.9]{dijk-koor:qsf}), the fact that $\ksi$ is a $\tau$-invariant
two-sided coideal implies that
$\Aksi$ is a $\ast$-subalgebra which is invariant under the right
$\Uq$-action on $\AU$ (or, equivalently,  the left coaction of
$\AU$ on itself). Important for the study of $\Aksi$ is the fact that $X=J^{\si}$
is a solution of the reflection equation
\begin{equation}\label{e:refl}
R_{12}X_1 R^{-1}_{12} X_2 = X_2 R_{21}^{-1}X_1 R_{21},
\end{equation}
where $R_{12}:=R$, $R_{21}:=PRP(=R^+)$, $X_1=X\ten \hbox{id}_V$ and
$X_2:=\hbox{id}_V\ten X$. This fact can be verified by direct computations. 
\begin{rem} 
The algebra $\Aksi$ can be considered as a quantized
coordinate ring on the complex Grassmannian $U(n)/(U(n-l)\times U(l))$ 
in the following way (see \cite{nds:gras} for more details).
The quantum space of $n\times n$ Hermitean matrices is 
defined as the algebra generated by $X=(x_{ij})_{ij}$ with relations given by the
reflection equation \eqref{e:refl}. It can be endowed with a $\ast$-structure and 
a left $\AU$-coaction (the quantum analogue of the adjoint action).
Since $J^{\si}$ is a solution of \eqref{e:refl} it gives rise to a ($\ast$-invariant) 
character of the quantized algebra of functions on the space of Hermitean matrices. 
In other words, $J^{\si}$ corresponds
to a classical point in the quantum space of Hermitean matrices.
Then $\Aksi$ may be considered as the quantized algebra of functions 
on the adjoint orbit of the classical point corresponding to
$J^{\si}$ (see \cite[Prop. 2.4]{nds:gras}). 
Since $J^{\si}$ has two different eigenvalues $1$ and $-q^{2\si}$ with
multiplicity $n-l$ and $l$ respectively, this quantum adjoint orbit is
isomorphic with the complex Grassmannian $U(n)/(U(n-l)\times U(l))$.
\end{rem}
The quantized function algebra $A_q(U/K)$ \eqref{AUKeerst} 
can formally be interpreted as the algebra 
$\Aksi$ with $\si\rightarrow \infty$. To make this a little bit more
explicit, we write
\begin{equation*}
L^+ = \begin{pmatrix} {}^{11}L^+ & {}^{12}L^+&  {}^{13}L^+ \\
                         0   &   {}^{22}L^+ &  {}^{23}L^+ \\
                        0   &   0   &        {}^{33}L^+
       \end{pmatrix}, \quad
L^- = \begin{pmatrix} {}^{11}L^- & 0 &  0  \\
                       {}^{21}L^-  &   {}^{22}L^- &  0 \\
                        {}^{31}L^- & {}^{32}L^-  &   {}^{33}L^-
      \end{pmatrix},
\end{equation*}
where ${}^{11}L^+$ is an $l\times l$ matrix, ${}^{22}L^+$ an
$(n-2l)\times (n-2l)$ matrix etc. Let $D$ be the
$l\times l$ matrix with $1$'s on the antidiagonal and $0$'s
everywhere else. The coefficients of the matrix
$L^+J^\si - J^\si L^-$ coincide with 
the coefficients of the following six matrices up to a sign:
\begin{align}\label{e:ksi-explicit}
 (i) \quad & q^\si(D\cdot {}^{31}L^- - {}^{13}L^+ \cdot D) +
     (1-q^{2\si})({}^{11}L^+ - {}^{11}L^-) \notag \\
(ii) \quad & {}^{12}L^+  + q^\si D\cdot{}^{32}L^-,\notag \\
(iii) \quad & {}^{23}L^+ \cdot q^{\si}D  + {}^{21}L^-,\notag \\
(iv) \quad & {}^{22}L^+ - {}^{22}L^-,\\
(v)  \quad & {}^{11}L^+\cdot q^{\si}D  - q^\si D\cdot {}^{33}L^-,\notag \\
(vi)  \quad & {}^{33}L^+\cdot q^{\si}D  - q^\si D\cdot {}^{11}L^-.\notag
\end{align}
Obviously, the coefficients of the following matrix are also contained in
$\gok^\si$:
\begin{equation}\label{e:ksi-explicit-b}
(vii) \quad  q^\si(D \cdot {}^{13}L^+  - {}^{31}L^- \cdot D)
      + (1-q^{2\si})({}^{33}L^+ - {}^{33}L^-).
\end{equation}
For later use, observe that the following elements 
of the ``Cartan subalgebra'' $U_q(\goh)$ belong to $\ksi$:
\begin{equation}\label{e:cartanelements}
L^+_{ii}-L^-_{ii}\; (l<i<l'),\quad
L^+_{ii}-L^-_{i'i'}\; (1\leq i\leq l),\quad 
L^-_{ii}-L^+_{i'i'}\; (1\leq i\leq l).
\end{equation}

It is clear from \eqref{e:ksi-explicit} and \eqref{e:ksi-explicit-b}
that, in the limit $\si\to\infty$, 
the matrices in (i)--(vii) tend either to zero or to the
following matrices
\begin{equation}\label{e:k-explicit}
{}^{11}L^+ - {}^{11}L^-, \quad {}^{12}L^+, \quad {}^{21}L^-, \quad
{}^{22}L^+ - {}^{22}L^-, \quad {}^{33}L^+ - {}^{33}L^-.
\end{equation}
Again, the subspace $\gok^\infty \subset \Uq$ spanned by
the coefficients of the matrices in \eqref{e:k-explicit} is a $\tau$-invariant
two-sided coideal.
Now, on the one hand, $\gok^\infty$-invariance in 
a left or right $\Uq$-module $M$ is 
obviously the same as invariance with respect 
to the Hopf $\ast$-sub\-al\-ge\-bra
\begin{equation*}
U_q(\gok):= 
U_q(\gog\gol(n-l,\CC))\ten U_q(\gog\gol(l,\CC))
\hookrightarrow U_q(\gog\gol(n,\CC)),
\end{equation*}
where invariance of $v\in M$ with respect to $u\in \Uq$ 
should be interpreted as $u \cdot v = \ep(u)\cdot v$
(if $M$ is a left $\Uq$-module). Using the Hopf algebra duality between $\AU$ and 
$\Uq$ it can be easily shown that invariance of $v\in M$ with respect to
$U_q(\gok)$ is the same as invariance with
respect to $\AK$ (cf.\ \cite[Prop. 1.12]{dijk-koor:qsf}). It follows
that $\AK$-invariance is equivalent to $\gok^\infty$-invariance, hence
$A_q(\gok^\infty\backslash U)= \AUK$. 
\begin{rem}\label{th:LJinf}
It should be observed that the matrix $J^\infty$ also
satisfies the reflection equation \eqref{e:refl}, but the subspace
spanned by the coefficients of the matrix $L^+J^\infty -
J^\infty L^-$ is strictly smaller than $\gok^\infty$ and of little use for the
purposes of this paper.
\end{rem}
The following lemma is now a direct consequence of the arguments given above.
\begin{lem}
Let $M$ be a finite-dimensional right $A_q(U)$-comodule with linear basis $\{v_i\}_i$.
Consider $M$ as left $\Uq$-module using the differentiated action
\eqref{diffactionX}. Suppose that $v_{\si}:=\sum_{i}c_i(\si)v_i$ 
\textup{(}$c_i(\si)\in {\mathbb{C}}$\textup{)} 
is a $\ksi$-fixed vector for all $\si\in {\mathbb{R}}$ 
and that $c_i:=\lim_{\si\rightarrow\infty}c_i(\si)$ exists for all $i$.
Then $\sum_ic_iv_i$ is an $A_q(K)$-fixed vector in $M$.
\end{lem}
\begin{rem}\label{th:cl-limit}
In some suitable algebraic sense (cf.\ 
\cite[Prop.\ 9.2.3]{cp:quantum}) the algebra $\Uq$ ``tends''
to $U(\gog)$ when $q$ tends to $1$. The corresponding 
limits of the $L$-o\-pe\-ra\-tors are given by
\begin{equation*}
L^\pm_{ij}/(q-q^{-1}) \to \pm e_{ji}\; (i \lessgtr j),\quad
(q^{\varepsilon_i}-q^{-\varepsilon_i})/(q-q^{-1}) \to e_{ii}
\end{equation*}
(cf.\ \cite[(1.10), (1.11)]{nou:macdonald}).
Hence, by \eqref{e:k-explicit} respectively \eqref{e:ksi-explicit},
the subspace $\gok^\si\subset \Uq$ ($\si=\infty$ respectively $\si=0$) tends
to the Lie subalgebra $\gok= \gog\gol({n-l},\CC)\oplus 
\gog\gol(l,\CC)\subset \gog$ respectively $\gok'\subset \gog$ (cf. \
\S 2) in the limit $q\to 1$. 
\end{rem}
Reflection equations play an important role
in the quantization of symmetric spaces 
(cf.\ \cite[\S 2]{nou:macdonald}, \cite{nou-sug:toyonaka}). 
For the purposes of this chapter, the importance of
this equation lies in the following fact. 
Recall that a vector $w$ in a left
$\Uq$-module $M$ is called $\ksi$-fixed if $\ksi\cdot w= 0$ (a similar
definition can be given for right $\Uq$-modules).

\begin{prop}\textup{(}\cite[Prop.\ 3.1]{NS1}, \cite{nds:gras}\textup{)} 
\label{th:wrefl}
Let $J$ be any $n\times n$ complex matrix. Write
$\gok^J\subset \Uq$  for the two-sided coideal spanned by the
coefficients of $L^+ J - J L^-$.
The element 
\begin{equation*}
w^J := \sum_{i,j} J_{ij}v_i\ten v_j^\ast\in V\ten V^\ast
\end{equation*}
in the left $\Uq$-module $V\ten V^\ast$ is a $\gok^J$-fixed vector if and
only if $J$ satisfies the reflection equation \textup{\eqref{e:refl}}.
\end{prop}
\begin{proof}
In the proof the same notational conventions as in 
\cite[Proof of Proposition 2.3]{nou:macdonald} will be used. 
Recall that the $\Uq$-module structure on 
$V^\ast$ corresponding to the dual $\AU$-comodule $V^\ast$ 
is given by
\begin{equation*}
u\cdot v^\ast(v) := v^\ast(S(u)\cdot v)\quad (u\in \Uq, v^\ast\in V^\ast,
v\in V). 
\end{equation*}
Set $\vek{v}:=(v_1, \ldots, v_n)$, then it follows from \eqref{e:Lvector} that
\begin{equation}\label{e:Lv}
L_1^\pm\cdot \vek{v}_2 = \vek{v}_2\cdot R^\pm_{12},\quad
L_1^+\cdot \vek{v}^\ast_2 = \vek{v}^\ast_2\cdot (R^-_{21})^{t_2},\quad
L_1^-\cdot \vek{v}^\ast_2 = \vek{v}^\ast_2\cdot (R_{21}^+)^{t_2}.
\end{equation}
Here ${}^{t_2}$ denotes transposition with respect to the second tensor
factor. An equation like  $L_1^+\cdot \vek{v}_2 = \vek{v}_2\cdot R^+_{12}$ 
should be interpreted as 
$L^+_{ij}\cdot v_k = \sum_{l=1}^n (R^+_{12})^{il}_{jk} v_l$ for all
$1\leq i,j,k\leq n$, where $R_{12}^+=\sum_{i,j,k,l}(R_{12}^+)^{ik}_{jl}e_{ij}\ten
e_{kl}$. Using the identities \eqref{e:Lv} one computes in shorthand notation,
\begin{equation*}
\begin{split}
L^+J\cdot w^J &= L^+_1 \cdot(\vek{v}_2 J_2\ten (\vek{v}^\ast)^t_2)J_1 =
(L_1^+\cdot \vek{v}_2) J_2\ten L_1^+\cdot (\vek{v}^\ast)^t_2J_1 \\
&= \vek{v}_2 R_{12}^+J_2R_{21}^-J_1 \ten (\vek{v}^\ast)^t_2,
\end{split}
\end{equation*}
since by \eqref{e:Lv} one has 
$L_1^+\cdot (\vek{v}^\ast)^t_2 = R_{21}^- \cdot (\vek{v}^\ast)^t_2$.
On the other hand,
\begin{equation*}
\begin{split}
JL^- \cdot w^J &= J_1L_1^- \cdot (\vek{v}_2J_2 \ten (\vek{v}^\ast)^t_2)=
 J_1 \vek{v}_2R_{12}^-J_2\ten L_1^-\cdot (\vek{v}^\ast)^t_2 \\
& = \vek{v}_2  J_1R_{12}^-J_2 R^+_{21} \ten (\vek{v}^\ast)^t_2,
\end{split}
\end{equation*}
since by \eqref{e:Lv} one has 
$L_1^-\cdot (\vek{v}^\ast)^t_2= R^+_{21}(\vek{v}^\ast)^t_2$.
It follows from the two preceding computations
that $w^J$ is $\gok^J$-fixed if and
only if $R_{12}^+J_2R_{21}^-J_1 = J_1 R_{12}^-J_2 R^+_{21}$.
Multiplying this last equation from the left and from the right 
by the permutation operator $P$ gives \eqref{e:refl}, which
proves the proposition.
\end{proof}
By Proposition \ref{th:wrefl} and the fact that the matrix $J^\si$ satisfies
the reflection equation \eqref{e:refl}, it follows that 
\begin{equation}\label{e:wsi}
w^\si:= \sum_{ij} J^\si_{ij}v_i\ten v_j^\ast\in V\ten V^\ast
\end{equation}
is a $\ksi$-fixed vector in the left $\Uq$-module $V\ten V^\ast$. Observe that 
$\lim_{\si\to\infty}w^\si =w^\infty$, with $w^\infty$ the right $\AUK$-fixed
vector defined by
\begin{equation}\label{e:winf}
w^{\infty}:=\sum_{ij}^{\infty}J_{ij}^{\infty}v_i\ten v_j^\ast=
\sum_{i=1}^{n-l}v_i\ten v_i^\ast.
\end{equation} 
Since $V\ten V^\ast\simeq V(\varpi_1)\oplus V(0)$ as left $\Uq$-modules (where $V(0)$
is the trivial module) and since $w^{\si}$ has a non-zero weight component 
of weight $\varpi_1$, it follows that $V(\varpi_1)$ 
has a non-zero $\ksi$-fixed vector. 

Next we construct a right $\ksi$-fixed vector in $V^\circ\ten (V^\ast)^\circ$.
Observe 
that a vector $\tilde{w}=\sum_{i,j}\tilde{J}_{ij}v_i\ten v_j^\ast\in V^\circ\ten 
(V^\ast)^\circ$ for a real matrix $\tilde{J}=\sum_{ij}\tilde{J}_{ij}e_{ij}$
is right $\ksi$-fixed if and only if $\tilde{w}$ is left $S(\ksi)$-fixed as element
in $V\ten V^\ast$ by the $\tau$-invariance of $\ksi$ and by Remark \ref{starleftright}. 
Reasoning as in the proof of Proposition \ref{th:wrefl},
it follows that $\tilde{w}$ is left $S(\ksi)$-fixed if $\tilde{J}$ is a solution
of the linear equation
\begin{equation}\label{refalt}
J_1^{\si}(R_{21}^-)^{t_1}\tilde{J}_2\bigl((R_{21}^-)^{t_1}\bigr)^{-1}=
R^{t_1}\tilde{J}_2\bigl(R^{t_1}\bigr)^{-1}J_1^{\si}
\end{equation}
where $J^{\si}$ is given by \eqref{Jsi}.
A solution $\tilde{J}=\tilde{J}^{\si}$ of \eqref{refalt} is given by
\begin{equation}\label{e:tildewsi}
\begin{split}
\tilde{J}^\sigma := \sum_{1\leq k\leq l} 
(1 - q^{2(n-2l)}q^{2\sigma})e_{kk} + &\sum_{l<k<l'} e_{kk}\\
-q^{\sigma-1}\sum_{1\leq k\leq l} &q^{2(k-l)}e_{kk'}
-q^{\sigma-1}\sum_{1\leq k\leq l} q^{2(k'-l)}e_{k'k}.
\end{split}
\end{equation}
We write $\tilde{w}^{\sigma}=\sum_{ij}\tilde{J}^{\si}_{ij}v_i\ten v_j^\ast$ for the
corresponding right $\ksi$-fixed vector in $V^\circ\ten (V^\ast)^\circ$.
In the same way as in the case of left $\ksi$-fixed vectors it follows that
$V(\varpi_1)^\circ$ has a non-zero right $\ksi$-fixed vector.
Observe that $\lim_{\si\rightarrow\infty}\tilde{w}^{\si}=w^{\infty}$, with
$w^{\infty}$ the $\AK$-fixed vector given by \eqref{e:winf}. 

Recall from the previous section that $V(\lambda)$ has at most one $\gok^\infty$-fixed
vector up to scalar multiples, and that $V(\lambda)$ has non-zero
$\gok^\infty$-fixed vectors if and only if $\lambda\in P_K^+$ 
(cf.\ Theorem \ref{th:gelfand}). We have the following analogous statement for
$\ksi$-fixed vectors $(-\infty<\si<\infty)$.
\begin{thm}\textup{(}\cite[Thm.\ 2.6]{nds:gras}\textup{)} \label{th:ksi-fixed}
Let $\lambda\in P^+$ and fix $-\infty<\si<\infty$. The irreducible 
left $\Uq$-module $V(\lambda)$ with highest weight $\lambda$ has
at most one $\ksi$-fixed vector \textup{(}up to scalar multiples\textup{)}.
There exist non-zero $\ksi$-vectors in $V(\lambda)$ 
if and only if $\lambda\in P^+_K$. The same statement holds for right $\ksi$-fixed
vectors in $V(\lambda)^\circ$.
\end{thm}
In the remainder of this section a proof of 
Theorem \ref{th:ksi-fixed} will be given.
Fix a parameter $-\infty < \si<\infty$.
First of all, we have the following crucial lemma. 
\begin{lem}\label{th:ksi-hwcomponent}
Let $\lambda\in P^+$ and fix $-\infty<\si<\infty$.
Then any non-zero $\ksi$-fixed vector
in the left $\Uq$-module $V(\lambda)$ 
has a non-zero weight component of highest weight $\lambda$.
The same statement holds for the right $\Uq$-module $V(\lambda)^{\circ}$.
\end{lem}
The proof of the lemma follows by analyzing the particular form 
of the two-sided coideal $\gok^\si$. The details are omitted here,
since the proof is analogous to the proof of 
\cite[Lemma 3.2]{nou:macdonald} and \cite[Prop. 3.2]{dijk-nou:proj}.

Since the vector subspace of $V(\lambda)$ (respectively
$V(\lambda)^\circ$) consisting of weight vectors
of weight $\lambda$ is one-dimensional, 
it follows from Lemma \ref{th:ksi-hwcomponent} that 
every irreducible finite-dimensional $P$-weighted $\Uq$-module has at most 
one $\ksi$-fixed vector up to scalar multiples.

Set $P_K=\oplus_{1\leq r\leq l}{\mathbb{Z}}\varpi_r$, where $\varpi_r$ are the
fundamental spherical weights (cf.\ \S 2).
Observe that the assignment $\lambda\mapsto \lambda^\natural$ as defined
in \S 2 extends to a bijection from
$P_K$ onto $P_{\Sigma}$. For $\mu\in P_{\Sigma}$, we write $\mu^\flat\in P_K$
for the inverse of $\mu$ under the bijection $\natural$. For later use let us record
the following elementary facts. Recall that $\FSW=\FSW_l$ denotes the Weyl
group of the root system $\Sigma$ (cf.\ \S 2).
\begin{lem}\label{weightorbits}
The bijection $\lambda\mapsto \lambda^\natural$ preserves
the dominance ordering on $P_K$ and $P_\Sigma$. If $\nu\in P_K$ then the image
under $\lambda\mapsto \lambda^\natural$ of $\bigl(\goS_n\cdot \nu\bigr) \cap P_K$
is equal to the $\FSW$-orbit $\FSW\cdot \nu^\natural$ in $P_\Sigma$.
\end{lem}
The following lemma is immediate from the fact that the Cartan type elements
listed in \eqref{e:cartanelements} belong to $\ksi$.
\begin{lem}\label{weightinPK}
Let $\lambda\in P_K^+$, $-\infty<\si<\infty$ and assume that 
$v\in V(\lambda)$ is a non-zero left $\ksi$-fixed
vector. Let $v=\sum_{\mu\leq\lambda}v_{\mu}$ be the decomposition of $v$
in weight vectors, where $v_{\mu}$ has weight $\mu\in P$. Then $v_{\mu}=0$
unless $\mu\in P_K$.
The same statement is valid for the right $\Uq$-module $V(\lambda)^{\circ}$.
\end{lem}
It follows from Lemma \ref{th:ksi-hwcomponent} and 
Lemma \ref{weightinPK} that if $V(\lambda)$ (respectively
$V(\lambda)^\circ$) has a non-zero $\ksi$-fixed vector, then $\lambda\in P_K^+$.

To finish the proof of Theorem \ref{th:ksi-fixed} we have to show that
all modules $V(\lambda)$ and $V(\lambda)^{\circ}$ 
($\lambda\in P_K^+$) have non-zero $\ksi$-fixed vectors. 
The existence of non-trivial $\ksi$-fixed vectors in
$V(\varpi_1)$ and in $V(\varpi_1)^{\circ}$ is already proved. 
Explicit intertwining operators
\begin{equation*}
\widehat{\Psi}_r\colon (V\ten V^\ast)^{\ten r} 
\to \Lambda_q^r(V)\ten \Lambda_q^r(V^\ast),\quad (1\leq r\leq l)
\end{equation*}
will be constructed to prove the existence of $\ksi$-fixed vectors
in higher fundamental spherical representations. The proof of Theorem
\ref{th:ksi-fixed} is then completed by computing the so-called principal term of 
$\widehat{\Psi}_r\bigl((w^{\sigma})^{\otimes r}\bigr)$, with $w^{\si}\in
V(\varpi_1)$ the $\ksi$-fixed vector given by \eqref{e:wsi}.

Before giving the construction of
$\widehat{\Psi}_r$, we first introduce the notion of principal term of a vector
$v\in\Lambda_q^r(V)\ten\Lambda_q^r(V^\ast)$ (cf.\ \cite{NS1}, \cite{S}).
For the present setting it is convenient to use a slightly modified definition of 
Noumi's and Sugitani's notion of principal term (cf.\ \cite{NS1}, \cite{S}). 
The definition is based on certain specific properties of the comodule
$\Lambda_q^r(V)\ten\Lambda_q^r(V^\ast)$. 
The comodule $\Lambda_q^r(V)\ten \Lambda_q^r(V^\ast)$ has a multiplicity-free 
decomposition
\begin{equation}\label{e:pieri}
\Lambda_q^r(V)\ten \Lambda_q^r(V^\ast) \cong \bigoplus_{s=0}^r 
V(\varpi_s)\quad (1\leq r\leq l)
\end{equation}
as right $A_q(U)$-comodules, where $\varpi_0:=0\in P_K^+$.
The decomposition \eqref{e:pieri} can be proved by computing the restriction
of the character of the module $\Lambda_q^r(V)\ten \Lambda_q^r(V^\ast)$ to the 
torus and using the classical Pieri formula for Schur functions 
\cite[I, (5.17)]{mac:symmetric} (cf.\ Proposition \ref{th:decomp}).
Due to the multiplicity-free decomposition \eqref{e:pieri}, the module
$\Lambda_q^r(V)\ten \Lambda_q^r(V^\ast)$ is very useful for the study of
$\ksi$-fixed vectors in $V(\varpi_r)$, as will be shown in the remainder of this
chapter as well as in the next chapter.
It follows from \eqref{e:pieri} that all the weights $\mu\in P$ of the module
$\Lambda_q^r(V)\ten \Lambda_q^r(V^\ast)$ are $\leq \varpi_r$, where
$\leq$ denotes the dominance order.
The vector $v_{[1,r]}\ten v_{[n-r+1,n]}^\ast\in\Lambda_q^r(V)\ten\Lambda_q^r(V^\ast)$
is the highest weight vector of the unique copy of $V(\varpi_r)$ within
$\Lambda_q^r(V)\ten\Lambda_q^r(V^\ast)$.
Suppose now that $v=\sum_{\mu\leq\varpi_r}v_{\mu}$ is the weight space decomposition of
a vector $v\in \Lambda_q^r(V)\ten\Lambda_q^r(V^\ast)$, 
where $v_{\mu}$ is the weight component of weight $\mu\in P$. 
Then the {\it principal term} of $v$ is defined by
\begin{equation}\label{leadingterm}
[v]:=\sum_{\nu\in \FSW(1^r)}v_{\nu^\flat}
\end{equation}
(cf.\ \cite{NS1}, \cite{S}), where $\FSW$ acts on $(1^r)\in P_{\Sigma}^+\subset
P_{\Sigma}={\mathbb{Z}}^l$ by permutations and sign changes (cf.\ \S 2).
It follows from Lemma \ref{weightorbits} that the principal term of a vector 
$v\in \Lambda_q^r(V)\ten\Lambda_q^r(V^\ast)$
lies in the unique copy of $V({\varpi_r})$ within
$\Lambda_q^r(V)\ten\Lambda_q^r(V^\ast)$.
By Lemma \ref{weightinPK} and Lemma \ref{weightorbits} one has:
\begin{lem}\label{convex}
Let $v$ be  a non-zero $\ksi$-fixed vector in $\Lambda_q^r(V)\ten\Lambda_q^r(V^\ast)$.
If $v-[v]$ has a non-zero weight component of weight $\nu$
then $\nu\in P_K$ and $\nu^\natural\in C(\varpi_r)$, where 
\begin{equation}\label{convexhull}
C(\mu):=\{ \mu'\in P_{\Sigma} \, | \, w\mu'< \mu \,\,\forall w\in \FSW \}\quad
(\mu\in P_{\Sigma}^+)
\end{equation}
is the strict integral convex hull of $\FSW\mu$.
\end{lem}
In the next proposition, the principal term of a $\ksi$-fixed vector in
$\Lambda_q^r(V)\ten\Lambda_q^r(V^\ast)$ (respectively in $\Lambda_q^r(V)^\circ
\ten\Lambda_q^r(V^\ast)^\circ$) is compared with the elements $u_r$, $\tilde{u}_r$ 
($1\leq r\leq l$) defined by
\begin{equation}\label{lr}
u_r:=\sum_{\stackrel{{\scriptstyle{I\subset [1,l]\cup [l',n]}}}
{{\scriptstyle{|I|=r, I\cap I'=\emptyset}}}}
v_I\ten v_{I'}^\ast,\quad 
\tilde{u}_r:=\sum_{\stackrel{{\scriptstyle{I\subset [1,l]\cup [l',n]}}}
{{\scriptstyle{|I|=r, I\cap I'=\emptyset}}}}
q^{\langle 2\rho, \tilde{\ep}_{I'}\rangle}v_I\ten v_{I'}^\ast
\end{equation}
where $I':=\{ i' \, | \, i\in I\}$. The element
$u_r$ lies in the unique copy 
of $V(\varpi_r)$ within $\Lambda_q^r(V)\ten\Lambda_q^l(V^\ast)$, whereas
$\tilde{u}_r$ lies in the unique copy of $V(\varpi_r)^\circ$
within $\Lambda_q^r(V)^\circ\ten\Lambda_q^r(V^\ast)^\circ$.
Observe that by the explicit form of the $\ksi$-fixed vectors $w^{\si}$ respectively
$\tilde{w}^{\si}$, we have
\begin{equation}\label{step1ind}
\lbrack w^{\si}\rbrack=-q^{\si}u_1,\quad 
\lbrack \tilde{w}^{\si}\rbrack=-q^{\si-1}q^{2(1-l)}\tilde{u}_1.
\end{equation}
For the construction of the intertwiner $\widehat{\Psi}_r$, consider now
the linear bijection $\beta\colon V^\ast \ten V\to V\ten V^\ast$ 
determined by
\begin{equation}\label{e:betadef}
\be(v^\ast_i\ten v_j)=q^{-\delta_{ij}}v_j\ten v_i^\ast +(q^{-1}-q)
\delta_{ij}\sum_{k<j} v_k\ten v^\ast_k.
\end{equation}
Write $V_i:=V$, $V^\ast_i:=V^\ast$ \textup{(}$1\leq i\leq r$\textup{)}.
Define a linear bijection
\begin{equation*}
\Psi_r\colon (V_1\ten V^\ast_1)\ten \cdots \ten (V_r\ten V^\ast_r)
\to (V_1\ten \cdots \ten V_r)\ten (V^\ast_1\ten\cdots\ten V^\ast_r)
\end{equation*}
by
\begin{equation}\label{e:Psidef}
\Psi_r := \be_{1,r} \circ \be_{2,r} \circ \cdots \circ \be_{r-1,r}
\circ \cdots \circ\be_{13}\circ \be_{23}\circ \be_{12},
\end{equation}
where $\beta_{ij}$  acts by definition as the identity
on all factors of the tensor product except for $V^\ast_i\ten V_j$, on which it 
is equal to $\beta$. Write
\begin{equation*}
\pr_r\colon V^{\ten r} \to \Lambda_q^r(V), \quad
\pr^\ast_r\colon (V^\ast)^{\ten r}\to \Lambda_q^r(V^\ast)
\end{equation*}
for the canonical projections. We now have the following generalization 
of \eqref{step1ind}.
\begin{prop}\label{th:intertwining}
Let $1\leq r\leq l$.
The operator 
\begin{equation*}
\widehat{\Psi}_r\colon 
(V\ten V^\ast)^{\ten r} \to \Lambda_q^r(V)\ten \Lambda_q^r(V^\ast)
\end{equation*}
defined by $\widehat{\Psi}_r := (\pr_r\ten \pr_r^\ast) \circ \Psi_r$
is a surjective intertwiner, and 
\begin{equation}
\begin{split}
\lbrack\widehat{\Psi}_r\bigl((w^\si)^{\otimes r}\bigr)\rbrack&=c_r({\si})u_r,\quad
c_r({\si}):=\left(\frac{q^{\si}}{q^2-1}\right)^r\bigl(q^2;q^2\bigr)_r,\nonumber\\
\lbrack \widehat{\Psi}_r\bigl((\tilde{w}^{\si})^{\otimes r}\bigr)\rbrack&=
\tilde{c}_r(\si)\tilde{u}_r,\quad 
\tilde{c}_r(\si):=\left(\frac{q^{\si-1}q^{2(1-l)}}{q^2-1}\right)^r\bigl(q^2;q^2\bigr)_r.
\nonumber
\end{split}
\end{equation}
\end{prop}
Before giving a proof of Proposition \ref{th:intertwining}, we first show how
it implies Theorem \ref{th:ksi-fixed}. 
Since $\ksi$ is a two-sided coideal and $\Psi_r$ an intertwining operator,
Proposition \ref{th:intertwining} shows that $\Psi_r\bigl((w^\si)^{\otimes r}\bigr)$
is a non-zero $\ksi$-fixed vector. 
Proposition \ref{th:intertwining} implies that the principal term of
$\Psi_r\bigl((w^\si)^{\otimes r}\bigr)$ is non-zero, hence 
it follows that $V(\varpi_r)$ ($1\leq r\leq l$) 
has a non-zero $\ksi$-fixed vector.
Since any $\lambda\in P_K^+$ can be written as a positive integral 
linear combination of the fundamental spherical weights 
$\{  \varpi_r\}_{1\leq r\leq l}$, it follows by an easy argument
using tensor products and Lemma \ref{th:ksi-hwcomponent}
that any $\lambda\in P^+_K$ is actually spherical. 
For right $\ksi$-fixed vectors the same argument holds,
since $\widehat{\Psi}_r$ is also an intertwiner as map from the module
$(V^\circ\ten (V^\ast)^\circ)^{\otimes r}$ to 
$\Lambda_q^r(V)^{\circ}\ten\Lambda_q^r(V^\ast)^{\circ}$ 
(cf.\ Remark \ref{th:Kleftright}).

So it remains to prove Proposition \ref{th:intertwining}.
The proof of this proposition, which proceeds by induction on $r$, 
is broken up into a couple of lemmas.
\begin{lem}\label{th:hatPhir}
For $2\leq r\leq n+1$ the linear mapping 
\begin{equation*}
\widehat{\Phi}_r\colon \Lambda_q^{r-1}(V^\ast)\ten V
\to V \ten \Lambda_q^{r-1}(V^\ast)
\end{equation*}
defined on the basis vectors
$v_I^\ast\ten v_j$ \textup{(}$|I|=r-1$, $1\leq j\leq n$\textup{)} by
\begin{equation*}
\widehat{\Phi}_r(v_I^\ast\ten v_j) =
\begin{cases}
v_j \ten v_I^\ast& \text{if $j\notin I$,}\\
q^{-1} v_j\ten v_I^\ast -(q-q^{-1}) {\displaystyle{\sum_{m<j}}}
{\displaystyle{\frac{\sgn(I\backslash j;m)}{\sgn(I\backslash j;j)}}} 
\,v_m\ten v^\ast_{(I\backslash j) \cup m}&\text{if $j\in I$}
\end{cases}
\end{equation*}
is an intertwining operator of right $A_q(U)$-comodules.
\end{lem}
\begin{proof}
Let $P\colon V\ten V\to V\ten V$ denote the flip. Define a linear
bijection $\gamma\colon V\ten V\to V\ten V$ by $\gamma:= PR$, 
with $R$ as in \eqref{e:Rdef}. 
The action of $\gamma$ on the basis vectors $v_i\ten v_j$ 
($1\leq i,j\leq n$) is given by
\begin{equation}\label{e:gammadef}
\gamma(v_i\ten v_j) = q^{\delta_{ij}}v_j\ten v_i + (q-q^{-1}) \theta_{i,j}
v_i\ten v_j
\end{equation}
with $\theta_{i,j}:=1$ if $i<j$ and $\theta_{i,j}:=0$ otherwise.
The fact that
the commutation relations between the $t_{ij}\in \AU$ can be written
as $RT_1T_2=T_2T_1R$ (cf.\ \S 3) implies
that $\gamma$ is an intertwining operator. 
Since $R$ is a solution of the Quantum Yang-Baxter Equation, $\gamma$ satisfies 
\begin{equation}\label{QYBE}
\gamma_{1}\circ \gamma_{2}\circ\gamma_{1}=
\gamma_{2}\circ\gamma_{1}\circ\gamma_{2},
\end{equation} 
with $\gamma_{i}\in \End (V^{\ten 3})$  acting
as $\gamma$ on the $i$th and $(i+1)$th tensor factors and as the identity on
the remaining factor.
Note furthermore that the exterior algebra $\Lambda_q(V)$ is isomorphic as a right
$A_q(U)$-comodule algebra with $T(V)/I$, where $T(V)$ is the tensor algebra of $V$
and $I\subset T(V)$ the two-sided ideal generated by 
$\ker(\id-q^{-1}\gamma)\subset V^{\otimes 2}\subset T(V)$.
Consider now the intertwiner $\Gamma_k\colon 
V^{\otimes (k-1)}\ten V\rightarrow V\ten \Lambda_q^{k-1}(V)$ 
($2\leq k\leq n+1$) defined by
\begin{equation*}
\Gamma_k=(\id\ten \pr_{k-1})\circ
\gamma_{1}\circ\gamma_{2}\circ\cdots\circ\gamma_{k-1}.
\end{equation*}
Application of \cite[Lemma 4.9 (1)]{has-hay} to
the Yang-Baxter operator $q^{-1}\gamma$ shows that there 
exists a unique bijective intertwiner 
\begin{equation*}
\widehat{\Gamma}_k\colon \Lambda_q^{k-1}(V)\ten V
\to V\ten \Lambda_q^{k-1}(V)
\end{equation*}
such that $\Gamma_k=\widehat{\Gamma}_k\circ (\pr_{k-1}\ten \id)$.
By a straightforward computation one verifies that
\begin{multline*}
\widehat{\Gamma}_k(v_I\ten v_j)= q^{|I\cap j|}v_j\ten v_I +\\
+ (1-q^2)(-q)^{-k+1}\sgn(I;j)
\sum_{\stackrel{{\scriptstyle{i\in I}}}{{\scriptstyle{i<j}}}}
\sgn(i;I\setminus i)v_i\ten v_{(I\setminus i)\cup j}
\end{multline*}
for $I\subset [1,n]$ with $|I|=k-1$ and $1\leq j\leq n$.

Next, the linear mapping 
$\delta_k: \Lambda_q^k(V^\ast) \to \Lambda_q^{n-k}(V)\ten \CC
\detq^{-1}$ ($1\leq k\leq n$) 
defined on the basis elements $v_I^\ast$ ($|I|=k$) by
$\delta_k(v_I^\ast):=\sgn(I;I^c)v_{I^c}\ten \detq^{-1}$
is a bijective intertwiner by \eqref{e:qminoranti}.
With the canonical identification $V\ten \CC\detq^{-1} \cong
\CC\detq^{-1} \ten V$ we have an intertwining operator
$\widehat{\Phi}_r\colon \Lambda_q^{r-1}(V^\ast)\ten V
\to V \ten \Lambda_q^{r-1}(V^\ast)$ defined by 
\begin{equation*}
\widehat{\Phi}_r := q^{-1}(\id\ten \delta_{r-1}^{-1})\circ
(\widehat{\Gamma}_{n-r+2}\ten \id)\circ (\delta_{r-1}\ten\id).
\end{equation*}
Starting from the explicit expressions for $\widehat{\Gamma}_{n-r+2}$ and
$\delta_{r-1}$, a straightforward calculation 
shows that $\widehat{\Phi}_r$ acts on the basis
vectors $v^\ast_I\ten v_j$ as required.
\end{proof}
\begin{cor}\label{th:corry}
The linear mappings $\beta$, $\Psi_r$, and $\widehat{\Psi}_r$
are right $A_q(U)$-comodule homomorphisms.
\end{cor}
\begin{proof}
The assertion follows from the previous lemma, since $\be=\widehat{\Phi}_2$
and the natural projections $\hbox{pr}_r$ and $\hbox{pr}_r^\ast$
intertwine the right $A_q(U)$-comodule actions.
\end{proof}
\begin{lem}\label{th:PhihatPhi}
Let $1\leq r\leq l$. The bijective intertwining operator
\begin{equation*}
\Phi_r\colon (V_1^\ast\otimes\cdots\otimes V_{r-1}^\ast)\ten V_1\to 
V_1\ten (V_1^\ast\otimes\cdots\otimes V_{r-1}^\ast)
\end{equation*}
defined by $\Phi_r := \be_{12} \circ \be_{23} \circ
\cdots \circ \be_{r-2,r-1}\circ \be_{r-1,r}$
satisfies 
\begin{equation*}
(\id\ten \pr_{r-1}^\ast)\circ \Phi_r = 
\widehat{\Phi}_r \circ (\pr^\ast_{r-1}\ten \id).
\end{equation*}
\end{lem}
\begin{proof}
For $I=\{i_1 <\ldots < i_r\}\subset [1,n]$, set
${\tilde{v}}_I^\ast:=
v_{i_r}^\ast\otimes\ldots\otimes v_{i_2}^\ast\otimes v_{i_1}^\ast$. 
It is clear from the definitions that
\begin{equation*}
(\id\ten \pr_{r-1}^\ast)\circ\Phi_r({\tilde{v}}_I^\ast 
\ten  v_j)
= v_j\ten v_I^\ast \quad \text{if $j\notin I$}.
\end{equation*}
If $j\in I$, then
\begin{equation*}
(\id\ten \pr_{r-1}^\ast)\circ\Phi_r({\tilde{v}}_I^\ast \ten v_j)
= q^{-1} v_j\ten v_I^\ast -(q-q^{-1}) \sum_{m<j} c(m,j) v_m\ten 
v^\ast_{(I\backslash j)\cup m},
\end{equation*}
where $c(m,j) :=(-q)^{|\{i\in I\mid m<i<j\}|}$ if $m\notin I$, and
$c(m,j):=0$ otherwise.
Using the definition of the $q$-signum $\sgn$, it follows that 
$c(m,j) = \sgn(I\backslash j;m)\sgn(I\backslash j; j)^{-1}$
if $m<j$, which concludes the proof of the lemma.
\end{proof}
Observe that the multiplication maps
\begin{equation*}
\mu\colon \Lambda_q(V)\ten\Lambda_q(V)\to \Lambda_q(V),
\quad \mu^*\colon \Lambda_q(V^\ast)\ten\Lambda_q(V^\ast)\to
\Lambda_q(V^\ast)
\end{equation*}
are intertwiners of the $\AU$-coactions, since
$\Lambda_q(V)$ and $\Lambda_q(V^\ast)$ are $A_q(U)$-co\-mo\-du\-le algebras.
\begin{lem}\label{th:hatThetar}
The intertwining operator
\begin{equation*}
\widehat{\Theta}_r \colon 
\Lambda_q^{r-1}(V) \ten \Lambda_q^{r-1}(V^\ast) \ten V\ten V^\ast
\to \Lambda_q^r(V)\ten \Lambda_q^r(V^\ast)
\end{equation*}
defined by $\widehat{\Theta}_r:= 
(\mu\ten \mu^\ast)\circ (\id_{\Lambda_q^{r-1}(V)}
\ten \widehat{\Phi}_r\ten \id_{V^\ast})$
satisfies 
\begin{equation}
\begin{split}
\lbrack\widehat{\Theta}_r(u_{r-1} \ten w^{\si})\rbrack &= 
-q^{\si}\frac{1-q^{2r}}{1-q^2}u_r\nonumber\\
\lbrack\widehat{\Theta}_r(\tilde{u}_{r-1}\ten \tilde{w}^{\si})\rbrack &=
-q^{\si-1}q^{2(1-l)}\frac{1-q^{2r}}{1-q^2}\tilde{u}_r\nonumber
\end{split}
\end{equation}
for $2\leq r\leq l$.
\end{lem}
\begin{proof}
If $v$ is a vector of weight $\mu$ in the domain of $\widehat{\Theta}_r$,
then $\widehat{\Theta}_r(v)$ is again a weight vector of weight $\mu$, since 
$\widehat{\Theta}_r$ intertwines the right $\AU$-coaction.
Hence, for a fixed
$I\subset [1,l]\cup [l',n]$ with $I\cap I'=\emptyset$ and $|I|=r-1$, we have
that $\lbrack\widehat{\Theta}_r(v_I\ten v_{I'}^\ast\ten v_s\ten v_t^\ast)\rbrack=0$
unless $s,t\not\in I\cup I'$ and $s\not=t$. By the explicit formulas for
the action of $\widehat{\Phi}_r$ (cf.\ Lemma \ref{th:hatPhir}), it follows that
\[
\lbrack\widehat{\Theta}_r(u_{r-1} \ten w^{\si})\rbrack=
-q^{\si}\sum_{I,k}v_I\wedge v_k\ten v_{I'}^\ast\wedge v_{k'}^\ast= 
-q^{\si}\sum_{J}
c_Jv_J\ten v_{J'}^\ast
\]
where the first sum is taken over pairs $(I,k)$ with 
$I\subset [1,l]\cup [l',n]$, $k\in [1,l]\cup [l',n]$, $|I|=r-1$, $I\cap I'=\emptyset$ 
and $k\not\in I\cup I'$, and the second sum is taken over subsets
$J\subset [1,l]\cup [l',n]$ with $J\cap J'=\emptyset$ and $|J|=r$.
The corresponding constant $c_J$ is given by
\[c_J=\sum_{k\in J}\sgn(J\setminus k;k)\sgn(k';J'\setminus k')=
\sum_{k\in J}\bigl(\sgn(J\setminus k;k)\bigr)^2=
\sum_{s=0}^{r-1}q^{2s}=\frac{1-q^{2r}}{1-q^2}.\]
The proof for the leading term of $\widehat{\Theta}_r(\tilde{u}_{r-1}\ten 
\tilde{w}^{\si})$ is similar.
\end{proof}
Proposition \ref{th:intertwining} can now be proved by induction to $r$,
using the previous lemma for the induction step.
\begin{proof}[Proof of Proposition \ref{th:intertwining}]
Define an intertwiner 
\begin{equation*}
\Theta_r\colon V^{\ten (r-1)}\ten (V^\ast)^{\ten(r-1)} \ten V\ten V^\ast\to
V^{\ten r}\ten (V^\ast)^{\ten r}
\end{equation*}
by 
\begin{equation*}
\Theta_r:= \id_{V^{\ten (r-1)}} \ten \Phi_r\ten \id_{V^\ast}.
\end{equation*}
It follows from Lemma \ref{th:PhihatPhi} that
\begin{equation}\label{e:ThetahatTheta}
(\pr_r\ten \pr^\ast_r)\circ \Theta_r =
\widehat{\Theta}_r \circ (\pr_{r-1}\ten \pr^\ast_{r-1}
\ten \id_V\ten \id_{V^\ast}).
\end{equation}
{}From the definitions of $\Psi_r$ and $\Phi_r$ it follows that
\begin{equation*}
\Psi_r = \Theta_r\circ (\Psi_{r-1}\ten\id)
\end{equation*}
and hence by \eqref{e:ThetahatTheta}
\begin{equation}\label{helftstuk}
\widehat{\Psi}_r = \widehat{\Theta}_r\circ
(\widehat{\Psi}_{r-1}\ten\hbox{id}).
\end{equation}
This allows us to prove the proposition by induction to $r$. The proposition
is trivial for $r=1$. Suppose that $r\geq 2$.
By the induction hypotheses and Lemma \ref{convex}
we have
\[\widehat{\Psi}_{r-1}\bigl((w^{\si})^{\ten r-1}\bigr)=c_{r-1}(\si)u_{r-1}+
\sum_{\nu\in C((1^{r-1}))}v_{\nu^\flat},\]
where $v_{\nu^\flat}$ is some weight vector of weight $\nu^\flat$ and 
$C(\mu)$ is defined by \eqref{convexhull}. For $\nu\in C((1^{r-1}))$ we have
$\lbrack \widehat{\Theta}_r(v_{\nu^\flat}\ten w^{\si})\rbrack=0$, hence
the induction step for the computation of $\lbrack
\widehat{\Psi}_r\bigl((w^{\si})^{\ten r}\bigr)\rbrack$ follows by combining 
Lemma \ref{th:hatThetar} with \eqref{helftstuk}. The leading term  
$\lbrack\widehat{\Psi}_r\bigl((\tilde{w}^{\si})^{\otimes r}\bigr)\rbrack$ can be
computed in a similar way.
\end{proof}
\begin{rem}\label{th:wrsi-hwcomponent}
It should be observed that the proof of Theorem \ref{th:ksi-fixed}
differs in important details from the proof of Theorem
\ref{th:gelfand}. Observe for instance that Lemma
\ref{th:ksi-hwcomponent} does not hold with $\ksi$-fixed
replaced by ${\mathfrak{k}}^{\infty}$-fixed, since any
${\mathfrak{k}}^{\infty}$-fixed
vector lies automatically in the zero weight space of the module. 
\end{rem}

%%%%%%%%%%%%%%%%%%%%%%%%%%%%%%%%%%%%%%%%%%%%%%%%%%%%%%%%%%%%%%%%%%%%%%%%%%%%%%%%%%
%%
%%
%%   SECTION 6   ZONAL SPHERICAL FUNCTIONS IN $\AUK$
%%
%%
%%
%%%%%%%%%%%%%%%%%%%%%%%%%%%%%%%%%%%%%%%%%%%%%%%%%%%%%%%%%%%%%%%%%%%%%%%%%%%%%%%%%%
\section{Zonal $(\si,\tau$)-spherical functions}
\label{section:zonal}
In this section the $\gok^\tau$-invariant 
($-\infty<\tau<\infty$) functions are studied in the
quantized coordinate ring $\Aksi$ $(-\infty<\si<\infty)$.
The results of this section were announced
in \cite[\S 3]{nds:gras}. 
The rank $1$ case of these results
were earlier derived by Koornwinder \cite{koor:su2} for $n=2$
and for arbitrary complex projective space by Noumi and Dijkhuizen
\cite{dijk-nou:proj}.

Let $-\infty < \si, \tau \leq \infty$ and denote $\Hst$ for the
$\ast$-subalgebra of left $\ksi$-invariant and right $\ktau$-invariant
functions in $\AU$. {}From Theorem \ref{th:gelfand}, 
Theorem \ref{th:ksi-fixed}
and \eqref{e:peterweyl} we obtain the decomposition
\begin{equation}
\Hst = \bigoplus_{\lambda\in P^+_K} \Hst(\lambda),
\quad \Hst(\lambda):= W(\lambda)\cap \Hst,
\end{equation}
the subspaces $\Hst(\lambda)$ ($\lambda\in P^+_K$) being one-dimensional.
A non-zero element $\varphi^{\si,\tau}(\lambda)\in \Hst(\lambda)$ is called a 
zonal $(\si,\tau)$-spherical function. Since the decomposition 
\eqref{e:peterweyl} is orthogonal with respect to
the inner product $\langle \varphi, \psi\rangle = h(\psi^\ast \varphi)$, 
the zonal spherical functions $\varphi^{\si,\tau}(\lambda)$ ($\lambda\in P^+_K$) 
are mutually orthogonal with respect to $\langle \cdot,\cdot \rangle$.
 
Let $M$ denote a right $\AU$-comodule with comodule mapping $\rho_M$ and
an invariant inner product $\langle \cdot,\cdot \rangle$. 
With any two elements $v,w\in M$ we associate the matrix coefficient
\begin{equation}\label{e:thetadef}
\theta_M(v,w) := \sum_{(w)} \langle w_{(1)}, v\rangle w_{(2)}\in \AU,
\quad \rho_M(w) =: \sum_{(w)} w_{(1)}\ten w_{(2)}.
\end{equation}
The map $\theta_M$ induces a linear map (denoted by the same symbol)
\begin{equation*}
\theta_M\colon M^\circ \ten M \to \AU,
\end{equation*}
which is surjective onto the subspace spanned by the matrix coefficients
of $M$. If no confusion can arise we sometimes write $\theta:=\theta_M$.
The following lemma is a direct consequence of 
these definitions (cf.\ \cite[Lemma 4.8]{nou:macdonald}).
\begin{lem}\label{th:theta}
Let $M$ be a unitary right $\AU$-comodule. The map $\theta_M\colon
M^\circ \ten M \to \AU$ satisfies the following properties:
\begin{enumerate}
\renewcommand{\labelenumi}{(\roman{enumi})}
\item $\theta_M$ is a $\AU$-bicomodule homomorphism, i.e.\ 
\begin{equation}
\Delta\circ \theta_M = (\theta_M\ten \id) \circ (\id\ten \rho_M), \quad 
\Delta\circ \theta_M = (\id\ten \theta_M) \circ (\rho_M^\circ \ten \id),
\end{equation}
where $\rho_M^\circ$ is defined as in Remark \ref{th:Kleftright}.
\item $\theta_M(v,w) = \tau(\theta_M(w,v))$ \textup{($v,w\in M$)}.
\item If $M$ is irreducible of highest weight
$\lambda\in P^+$, then $\theta_M\colon M^\circ\ten M \to W(\lambda)$ is an
isomorphism of $\AU$-bicomodules.
\end{enumerate}
\end{lem}
Lemma \ref{th:theta} can be used to construct zonal $(\si,\tau)$-spherical functions
as follows.
Let $v_{\si}(\lambda)\in V(\lambda)$ respectively 
$\tilde{v}_{\tau}(\lambda)\in V(\lambda)^\circ$ be a non-zero $\ksi$-fixed
respectively $\ktau$-fixed vector ($\lambda\in P_K^+$).  
Let $\langle \cdot,\cdot \rangle$ be a positive definite inner product
on $V(\lambda)$, and write $\theta_{\lambda}$ for the map $\theta$ in Lemma
\ref{th:theta} with respect to the unitary comodule $(V(\lambda),\langle \cdot,\cdot
\rangle)$. Then 
\begin{equation}\label{spher}
\varphi^{\sigma,\tau}(\lambda):=\theta_{\lambda}(\tilde{v}_{\tau}(\lambda),
v_{\si}(\lambda))\in
\Hst(\lambda)
\end{equation}
is a zonal $(\si,\tau)$-spherical function by Lemma \ref{th:theta}.
This leads to the following lemma.
\begin{lem}\label{landsinx}
Let $-\infty<\si,\tau<\infty$ and $\lambda\in P_K^+$.
The image of $\varphi^{\si,\tau}(\lambda)$ under
the restriction map ${}_{|\TT}\colon \AU\to \AT$ is of the form
\begin{equation}\label{specform}
{\varphi^{\si,\tau}(\lambda)}_{|\TT}=c_{\lambda^\natural}m_{\lambda^\natural}(x)
+\sum_{\nu\in C(\lambda^\natural)}c_{\nu}x^{\nu},\quad c_{\nu}\in {\mathbb{C}},
\end{equation}
with $c_{\lambda^\natural}\not=0$ and $C(\nu)$ given by \eqref{convexhull}.
Here the notation $x^{\nu}:=x_1^{\nu_1}x_2^{\nu_2}\ldots x_l^{\nu_l}$ for
$\nu=(\nu_1,\ldots,\nu_l)\in P_{\Sigma}$ is used, with
the $x_i$ \textup{($1\leq i\leq l$)} being defined by \eqref{e:xdef}.
\end{lem}
\begin{proof}
Since any $\lambda\in P_K^+$ can be written as a positive integral 
linear combination of the fundamental spherical weights 
$\{  \varpi_r\}_{1\leq r\leq l}$, it follows by an easy argument
using tensor products and Lemma \ref{th:ksi-hwcomponent}
that \eqref{specform} for arbitrary $\lambda\in P_K^+$ 
follows from \eqref{specform} for the fundamental
spherical weights $\lbrace \varpi_r\rbrace_{r=1}^l$. 

So fix a fundamental weight $\varpi_r$ ($1\leq r\leq l$). 
Consider the $\AU$-invariant inner product 
\begin{equation}\label{inner}
\langle v_I\ten v_J^\ast, v_K\ten v_L^\ast\rangle=
q^{-\langle 2\rho,\tilde{\ep}_J\rangle}\delta_{I,K}\delta_{J,L}
\end{equation}
on $\Lambda_q^r(V)\ten\Lambda_q^r(V^\ast)$ (cf. \ \S 3) 
and write $\theta$ for the map \eqref{e:thetadef} associated with the unitary comodule
$\bigl(\Lambda_q^r(V)\ten\Lambda_q^r(V^\ast),\langle\cdot,\cdot \rangle\bigr)$.
By \eqref{e:pieri}, the comodule $V(\varpi_r)$ may be considered as an irreducible
component of $\Lambda_q^r(V)\ten\Lambda_q^r(V^\ast)$ with invariant inner product
given by the restriction of $\langle\cdot,\cdot\rangle$ to $V(\varpi_r)$. 
Then, by Proposition \ref{th:intertwining}
and the fact that $u_r\in\Lambda_q^r(V)\ten\Lambda_q^r(V^\ast)$ 
(respectively $\tilde{u}_r\in\Lambda_q^r(V)^\circ\ten
\Lambda_q^r(V^\ast)^\circ$) lies in the unique irreducible component
$V(\varpi_r)$ (respectively $V(\varpi_r)^\circ$), the principal terms of the
$\ksi$-fixed vector $v_{\si}(\varpi_r)$ and the $\ktau$-fixed vector
$\tilde{v}_{\tau}(\varpi_r)$ are given by
\begin{equation}\label{hrest}
\lbrack v_{\si}(\varpi_r)\rbrack=c_ru_r,\quad
\lbrack\tilde{v}_{\tau}(\varpi_r)\rbrack=\tilde{c}_r\tilde{u}_r
\end{equation}
for non-zero constants $c_r,\tilde{c}_r\in \CC$.
For $v_{\mu}\in\Lambda_q^r(V)\ten\Lambda_q^r(V^\ast)$ of weight $\mu$
and $\tilde{v}_{\nu}\in\Lambda_q^r(V)^\circ\ten\Lambda_q^r(V^\ast)^\circ$ of
weight $\nu$ we have ${\theta(\tilde{v}_{\nu},v_{\mu})}_{|\TT}=0$ 
if $\mu\not=\nu$, and ${\theta(\tilde{v}_{\mu},v_{\mu})}_{|\TT}$ 
is a multiple of $z^{\mu}$.
Using Lemma \ref{weightinPK} and 
the fact that ${\mathbb{C}}[x^{\pm 1}]$
is the subalgebra of $A(\TT)$ spanned by the monomials $z^{\mu}(=x^{\mu^\natural})$ 
($\mu\in P_K$), we obtain from \eqref{hrest} that
\begin{equation}
\begin{split}
{\varphi^{\si,\tau}(\lambda)}_{|\TT}
&={\theta\bigl(\tilde{v}_{\tau}(\varpi_r),v_{\si}(\varpi_r)\bigr)}_{|\TT}\nonumber\\
&= {\theta\bigl(\lbrack\tilde{v}_{\tau}(\varpi_r)\rbrack, \lbrack v_{\si}(\varpi_r)
\rbrack\bigr)}_{|\TT}+\sum_{\nu\in C((1^r))}d_{\nu}x^{\nu}\nonumber\\
&= d_{(1^r)}m_{(1^r)}(x)+\sum_{\nu\in C((1^r))}d_{\nu}x^{\nu}\nonumber 
\end{split}
\end{equation}
with $d_{(1^r)}=c_r\tilde{c}_r\not=0$, since 
${\theta(\tilde{u}_r,u_r)}_{|\TT}=m_{(1^r)}(x)$.
This completes the proof of \eqref{specform} for the fundamental spherical
weights. 
\end{proof}
Lemma \ref{landsinx} has the following important consequence.
\begin{cor}\label{commuterende}
The restriction of the map ${}_{|\TT}\colon\AU\to A(\TT)$ to $\Hst$ 
defines an injection from $\Hst$ into
${\mathbb{C}}[x^{\pm 1}]$ for $-\infty<\si,\tau<\infty$. 
In particular, $\Hst$ is a commutative algebra for $-\infty<\si,\tau<\infty$.
\end{cor}
Recall from \cite{rtf:R} and \cite[\S 3]{nds:gras} the Casimir operator
\begin{equation*}
C := \sum_{ij} q^{2(n-i)}L^+_{ij}S(L^-_{ji})\in \Uq.
\end{equation*}
Since $C$ is central, it acts on $W(\lambda)$ ($\lambda\in P^+$) as
a scalar $\chi_\lambda(C)$, which is given by
\begin{equation*}
\chi_\lambda = \sum_{k=1}^n q^{2(\lambda_k+n-k)}.
\end{equation*}
Also, $C$ maps $\Hst$ into itself. Therefore, if $-\infty <\sigma,\tau<\infty$, 
the restricted Casimir operator $C\colon \Hst\to \Hst$ induces an operator 
\begin{equation*}
L\colon {\Hst}_{|\TT} \to {\Hst}_{|\TT}\subset \CC[x^{\pm 1}],
\end{equation*}
which is called the radial part of $C$. 
Explicitly, $L$ is the map satisfying
\[
L({\varphi}_{|\TT})={(C\varphi)}_{|\TT},\quad \forall \varphi\in\Hst.
\]
Crucial for the identification of the zonal $(\si,\tau)$-spherical functions 
is the identification of the radial part $L$ of the Casimir element $C$ 
with the restriction to ${\Hst}_{|\TT}$
of an explicit second-order $q^2$-difference operator on $\CC[x^{\pm 1}]$.  
Without proof we will state here the result (see \cite[\S 3]{nds:gras}).
\begin{thm}\textup{(}\cite{nds:gras}\textup{)}\label{th:radialpart} 
Let $-\infty < \si, \tau < \infty$ and $\lambda\in P^+_K$.
The operator $L-\chi_\lambda(C)\id$ coincides on ${\Hst}_{|\TT}\subset \CC[x^{\pm 1}]$ 
with a constant multiple of Koornwinder's second-order $q^2$-difference operator
$D - E_{\lambda^\natural}\id$ in the variables $x=(x_1,\ldots,x_l)$ 
with base $q^2$ and parameters $(\underline{t},t)=(\underline{t}^{\si,\tau},q^2)$,
given by  
\begin{equation}\label{e:parameters}
\begin{split}
t_0^{\si,\tau} &= -q^{\si + \tau + 1},\quad t_1^{\si,\tau} = -q^{-\si - \tau + 1},\\
t_2^{\si,\tau} &= q^{\si - \tau + 1}, \,\,\quad 
t_3^{\si,\tau} = q^{-\si + \tau + 2(n-2l) + 1}.
\end{split}
\end{equation}
\end{thm}
For a proof of the theorem for rank 1, see \cite{dijk-nou:proj}. 
In \cite{nou-sug:toyonaka} a proof can be found for the special case
$n=2l$ and $\si=\tau=0$.

Observe that for $-\infty<\si,\tau<\infty$ we have 
$\underline{t}^{\si,\tau} \in V_{K}$ (cf. \S \ref{section:polynomials}) and 
$t_0^{\si,\tau}t_1^{\si,\tau}t_2^{\si,\tau}t_3^{\si,\tau}\in (0,1)$. 
In particular, the eigenvalues $E_{\lambda^\natural}$ are 
mutually different for
compatible weights when $-\infty<\si,\tau<\infty$ (see \S \ref{section:polynomials}).

We write $D_{\si,\tau}$ for Koornwinder's 
second-order $q^2$-difference operator in base $q^2$ with parameters
$(\underline{t},t)=(\underline{t}^{\si,\tau},q^2)$, and we write $E_{\mu}^{\si,\tau}$ 
$(\mu\in P_{\Sigma}^+)$ for the
corresponding eigenvalues. We furthermore write
$P_{\mu}^{\si,\tau}(x):=P_{\mu}^K(x;\underline{t}^{\si,\tau};q^2,q^2)$ 
($\mu\in P_{\Sigma}^+$) for the corresponding monic Koornwinder polynomials. 
By Theorem \ref{th:radialpart},
${\varphi^{\si,\tau}(\lambda)}_{|\TT}\in {\mathbb{C}}[x^{\pm 1}]$ 
is an eigenfunction of $D_{\si,\tau}$ with eigenvalue
$E_{\lambda^\natural}^{\si,\tau}$ for $\lambda\in P_K^+$.
By \cite[Lemma 6.2]{nou-sug:toyonaka}, any  
eigenfunction $\varphi(x)\in {\mathbb{C}}[x^{\pm 1}]$
of $D_{\si,\tau}$ with eigenvalue $E_{\mu}^{\si,\tau}$ ($\mu\in P_{\Sigma}^+$)
and which is of the particular form
\[\varphi(x)=c_{\mu}m_{\mu}(x)+\sum_{\nu\in C(\mu)}c_{\nu}x^{\nu}\]
is a constant multiple of the Koornwinder polynomial $P_{\mu}^{\si,\tau}(x)$.
Combined with Lemma \ref{landsinx}, the following main result 
of the paper \cite{nds:gras} is obtained.
\begin{thm}\textup{(}\cite{nds:gras}\textup{)}\label{th:AWspher}
Let $-\infty < \si, \tau < \infty$. The restriction 
$\varphi^{\si,\tau}(\lambda)_{|\TT}$ of the zonal spherical function
$\varphi^{\si,\tau}(\lambda)\in \Hst(\lambda)$ \textup{(}$\lambda\in P^+_K$\textup{)} 
is equal to the Koornwinder polynomial
$P_{\lambda^\natural}^{\si,\tau}(x)$, up to a non-zero scalar multiple.
In particular, ${}_{|\TT}$ defines an algebra isomorphism from $\Hst$ onto
${\mathbb{C}}[x^{\pm 1}]^{\FSW}$.
\end{thm}

\begin{rem}\label{th:finitesigma}
It should be observed here that the assumption $-\infty < \si,\tau < \infty$
in the preceding arguments is absolutely essential. In fact, the map
${}_{|\TT}\colon \AU \to \AT$ factors through the projection
$\pi_K\colon \AU \to \AK$. This implies that the image
of $\Hst$ under ${}_{|\TT}$ is one-dimensional as soon as either 
$\si$ or $\tau$ is infinite.
\end{rem}
%%%%%%%%%%%%%%%%%%%%%%%%%%%%%%%%%%%%%%%%%%%%%%%%%%%%%%%%%%%%%%%%%%%%%%%%%%%%%%%%
%%
%%
%%            Section: Limit transitions on the quantum Grassmannian
%%
%%
%%%%%%%%%%%%%%%%%%%%%%%%%%%%%%%%%%%%%%%%%%%%%%%%%%%%%%%%%%%%%%%%%%%%%%%%%%%%%%%%%

\section{Limit transitions on quantum Grassmannians}
\label{section:limGras}
In this section we study the right $\gok^\tau$-invariant 
($-\infty<\tau\leq\infty$) functions in the
quantized coordinate ring $\AUK=A_q(\gok^\infty\backslash U)$.
Our method will be to regard this case as a limit of the case
$-\infty<\si,\tau<\infty$ by sending $\si$ to infinity. 
This limit can be made rigorous by 
using explicit information about the limit transitions from 
Koornwinder polynomials to multivariable big and little $q$-Jacobi polynomials. 
In the rank $1$ case these results 
were derived earlier by Koornwinder \cite{koor:su2} for $2$-spheres
and for arbitrary complex projective space by Noumi and Dijkhuizen
\cite{dijk-nou:proj}.

For the proper interpretation of the limit transitions of the zonal spherical
functions, a careful study is needed of 
the pre-images of the $\FSW$-invariant functions
$e_s(x):=m_{(1^s)}(x)$ ($1\leq s\leq l$) under the isomorphism
${}_{|\TT}: \Hst\to \CC[x^{\pm 1}]^{\FSW}$. 
For $-\infty < \si,\tau <\infty$, write $e_r^{\si,\tau}$ 
for the unique element in $\Hst$ such that its restriction to the torus
is equal to $e_r(x)$ ($1\leq r\leq l$). 
It is convenient to put $e_0(x):=1$ and $e_0^{\si,\tau}:=1$. 
Recall that the $\FSW$-invariant 
functions $\lbrace e_r\rbrace_{r=1}^l$ are algebraically
independent generators of $\CC[x^{\pm 1}]^{\FSW}$ 
(cf.\ \S\ref{section:polynomials}). 
In other words, by setting
\begin{equation}\label{Wgeval}
\hat{P}(e_1(x),\ldots,e_l(x)):=P(x),
\quad P\in  {\mathbb{C}}[x^{\pm 1}]^{\FSW}
\end{equation}
we get an algebra isomorphism $P\mapsto \hat{P}$ of ${\mathbb{C}}[x^{\pm 1}]^{\FSW}$
onto ${\mathbb{C}}[y]$, where $y=(y_1,\ldots,y_l)$ is an $l$-tuple of independent
variables. It follows from Theorem \ref{th:AWspher} that the elements
$\lbrace {e_r^{\si,\tau}}\rbrace_{r=1}^l$ are algebraically independent
generators of the algebra $\Hst$. 

Using Theorem \ref{th:AWspher} it is now easy to derive an explicit
form of the restriction of the normalized Haar functional $h$ to $\Hst$.
Recall that the parameters $\underline{t}^{\si,\tau}$ lie in the
parameter domain $V_{K}$ for $-\infty<\si,\tau<\infty$ (see \S
\ref{section:polynomials} for the definition of $V_{K}$).
Write $\langle \varphi\rangle_{\si,\tau}:=\langle \varphi,1 
\rangle^{\underline{t}^{\si,\tau}}_{K,q^2,q^2}$ for the constant term of 
$\varphi\in {\mathbb{C}}[x^{\pm 1}]^{\FSW}$, with $\langle\cdot,\cdot \rangle_K$
defined by \eqref{e:AWorth}.
Observe that $\langle 1\rangle_{\si,\tau}$ is non-zero by the positive
definiteness of $\langle\cdot,\cdot \rangle_K$. Explictly, $\langle 1\rangle_{\si,\tau}$
can be given explicitly as product and quotient of $q$-Gamma functions 
by Gustafson's evaluation of the multidimensional Askey-Wilson integral 
\eqref{e:AWintegral}. 
\begin{cor}\label{haarAW} Let $-\infty<\si,\tau<\infty$.
The restricted Haar functional $h: \Hst\rightarrow \CC$ is explicitly given
by
\[h\bigl(\hat{P}(e_1^{\si,\tau},e_2^{\si,\tau},\ldots,e_l^{\si,\tau})\bigr)=
\frac{\langle P \rangle_{\si,\tau}}{\langle 1\rangle_{\si,\tau}}
\quad (P\in \CC[x^{\pm 1}]^{\FSW}).
\]
\end{cor}
\begin{proof}
The left- and the right-hand side are both equal to zero for 
$P=P_{\mu}^{\si,\tau}$ with $0\not=\mu\in P_{\Sigma}^+$, and 
both equal to $1$ for $P=1$. The corollary follows now by linearity, 
since the Koornwinder polynomials $P_{\mu}^{\si,\tau}(x)$ ($\mu\in P_{\Sigma}^+$)
form a linear basis of $\CC[x^{\pm 1}]^{\FSW}$. 
\end{proof}
Recall the intertwiner 
\begin{equation*}
\widehat{\Psi}_r\colon 
(V\ten V^\ast)^{\ten r} \to \Lambda_q^r(V)\ten \Lambda_q^r(V^\ast)
\end{equation*}
defined in Proposition \ref{th:intertwining}. Introduce 
left $\gok^\si$-fixed vectors $w_r^\si\in \Lambda_q^r(V)\ten \Lambda_q^r(V^\ast)$
and right $\gok^\tau$-fixed vectors 
$\tilde{w}_r^\tau\in \Lambda_q^r(V)^\circ\ten \Lambda_q^r(V^\ast)^\circ$ by
\begin{equation*}
w_r^\si:= \widehat{\Psi}_r\bigl((w^\si)^{\otimes r}\bigr),\quad
\tilde{w}_r^\tau := \widehat{\Psi}_r\bigl((\tilde{w}^\tau)^{\otimes r}\bigr) 
\quad (1\leq r\leq l).
\end{equation*}
Here we have used the notation $\tilde{w}^{\infty}:=w^{\infty}$ \eqref{e:winf}
when $\tau=\infty$, which is consistent with the definition of $\tilde{w}^{\tau}$
for $-\infty<\tau<\infty$ since 
$\lim_{\tau\rightarrow\infty}\tilde{w}^{\tau}=w^{\infty}$.
Consider now the $(\si,\tau)$-spherical elements
\begin{equation}\label{e:philambda-xexpansion}
\varphi_r^{\si,\tau} := \theta(\tilde{w}_r^\tau, w_r^\si)\in
\bigoplus_{s=0}^r\Hst(\varpi_s)
\quad (1\leq r\leq l)
\end{equation}
(cf.\ \eqref{e:pieri}), where $\theta$ is the map \eqref{e:thetadef}
associated with the unitary comodule
$\Lambda_q^r(V)\ten \Lambda_q^r(V^\ast)$ endowed with the inner
product $\langle\cdot,\cdot \rangle$ (see \eqref{inner}
for the definition of $\langle \cdot,\cdot \rangle$).
It is convenient to put $\varphi_0^{\si,\tau} :=1$.
By Theorem \ref{th:AWspher} and \eqref{e:philambda-xexpansion},  
${\varphi_r^{\si,\tau}}_{|\TT}$ is a linear combination of
the $\FSW$-invariant functions $e_s(x)\in \CC[x^{\pm 1}]^{\FSW}$ ($0\leq s\leq r$).
\begin{lem}\label{th:limitlemma2}
Let $-\infty < \si,\tau < \infty$, $1\leq r\leq l$.
In the expansion
\begin{equation*}
{\phist}_{|\TT} = a^r_r(q^\sigma, q^\tau) e_r + \cdots +
a^r_0(q^\sigma, q^\tau) e_0,
\end{equation*}
each coefficient $a_i^r$ is a polynomial in $q^\si$ and $q^\tau$
which is the sum of monomials of partial degree $\geq i$ in each
of the variables. Moreover, $a_r^r(q^\sigma, q^\tau)= 
cq^{r\sigma+r\tau}$ with $c\neq 0$ independent of $q^{\si}$ and $q^{\tau}$.
\end{lem}
\begin{proof}
It is obvious from the definitions that
the coefficients are polynomial in $q^\si$ and $q^\tau$.
To prove the estimates on the partial degrees, we study
the action of the intertwiner $\widehat{\Psi}_r$ on the vectors
$(w^\si)^{\ten r}$ and $(\tilde{w}^\tau)^{\ten r}$ in detail. 
We proceed in a number of steps.

1) Let $1\leq i_1\leq \cdots \leq i_r\leq n$ and 
$1\leq j_1\leq \cdots \leq j_r\leq n$ be integers. We use 
the shorthand notation
$\underline{i}:= (i_1, \ldots, i_r)$, $\underline{j} :=
(j_1, \ldots, j_r)$. Call a tensor $t$ in some tensor
product space made up of factors $V$ or $V^\ast$ (the total
number of factors $V$ being equal to the total number of
factors $V^\ast$) a basic tensor of type $(\underline{i}, 
\underline{j})$ if $t$ is the tensor product in any given order
of the vectors
$v_{i_1}, \ldots, v_{i_r}$ and $v^\ast_{j_1}, \ldots,
v^\ast_{j_r}$. Let 
$n_k(\underline{i})$ denote the cardinality of the
set $\{p\in [1,r] \mid i_p=k\}$. For a
basic tensor $t$ of type $(\underline{i},\underline{j})$
define
\begin{equation*}
n(t):= \sum_{k=1}^n \min(n_k(\underline{i}), n_k(\underline{j})).
\end{equation*}
{}From an informal point of view, 
$n(t)$ is the number of factors $v_i$ in $t$
that ``cancel'' against a factor $v_i^\ast$. Recall the intertwiner
$\Psi_r\colon (V\ten V^\ast)^{\ten r} \to V^{\ten r}\ten (V^\ast)^{\ten r}$
defined in \eqref{e:Psidef}. Let $t$ be a basic tensor in 
$(V\ten V^\ast)^{\ten r}$. 
Since $\Psi_r$ is a composition of intertwiners $\beta_{ij}$ (see \eqref{e:Psidef})
it follows by inspection of \eqref{e:betadef} that $\Psi_r(t)$ is a linear combination
of basic tensors $t'$ in $V^{\ten r}\ten (V^\ast)^{\ten r}$ with $n(t')=n(t)$.

(2) A basic tensor $t\in (V\ten V^\ast)^{\ten r}$ is called typical
if it is a product of tensors in $V\ten V^\ast$ of type 
$v_i\ten v_i^\ast$ ($1\leq i\leq n-l$), $v_i\ten v_{i'}^\ast$,
$v_{i'}\ten v_i^\ast$ ($1\leq i\leq l$). We call a typical
tensor $t\in (V\ten V^\ast)^{\ten r}$ $k$-typical if the
number of factors of type $v_i\ten v_{i'}^\ast$ ($i\in [1,l]\cup [l',n]$) 
is equal to $k$. If $t$ is a $k$-typical tensor then $\widehat{\Psi}_r(t)$
is a linear combination of elements $v_I\ten v_J^\ast$ where
$I,J\subset [1,n]$ are such that $|I|=|J| =r$ and
$|I\cap J| \geq r-k$.
In fact, this follows from (1) and the definition of 
$\widehat{\Psi}_r$, since $n(t)\geq r-k$.

(3) It is an immediate consequence of 
the definition of the coactions on $\Lambda_q^r(V)$ and 
$\Lambda_q^r(V^\ast)$ and of 
\eqref{e:qminoranti} that
\begin{equation*}
\theta(v_I\ten v^\ast_J, v_K\ten v_L^\ast)_{|\TT} =
q^{-\langle 2\rho, \tilde{\ep}_J\rangle} \delta_{I,K} \delta_{J,L}
z^{\tilde{\ep}_I-\tilde{\ep}_J},
\end{equation*}
for $I,J\subset [1,n]$ with $|I|=|J|=r$.

(4)  Let $t$ be a $k$-typical tensor and $t'$ a $m$-typical tensor.
Let $\mu\in \FSW (1^i)\subset P_{\Sigma}$ ($i\in [1,r]$) be any weight,
and suppose that the coefficient of $z^{\mu^\flat}$ in the expansion of 
${\theta\bigl(\widehat{\Psi}_r(t), \widehat{\Psi}_r(t')\bigr)}_{|\TT}$
with respect to the basis $\lbrace z^{\lambda}\rbrace_{\lambda\in P}$ of $A(\TT)$,
is non-zero. Then $k\geq i$ and $m\geq i$. This is a 
straightforward consequence of (2) and (3).

(5) There is a unique expansion
$(w^{\sigma})^{\otimes r}=\sum_{k=0}^r\sum_{t_k}c_{t_k}t_k$
where $t_k$ runs over all $k$-typical tensors in $(V\ten V^\ast)^{\otimes r}$.
The non-zero  $c_{t_k}$ are linear combinations  of monomials $(q^{\si})^i$
with $i\geq k$. Similarly, there is a unique expansion
$(\tilde{w}^{\tau})^{\otimes r}=\sum_{k=0}^r\sum_{t_k'}d_{t_k'}t_k'$ 
where $t_k$ runs over all 
$k$-typical tensors in $(V\ten V^\ast)^{\otimes r}$.
The non-zero $d_{t_k'}$ are linear combinations of monomials $(q^{\tau})^i$
with $i\geq k$. Hence,
\begin{equation*}
{\phist}_{|\TT} = \sum_{k,m=0}^r\sum_{t_k,t_m'} c_{t_kt_m'} 
\theta(\widehat{\Psi}_r(t_m'),\widehat{\Psi}_r(t_k))_{|\TT}
\end{equation*}
with $c_{t_kt_m'}=c_{t_k}\overline{d_{t_m'}}$ 
a linear combination of monomials
$(q^{\si})^i(q^{\tau})^j$ with $i\geq k$ and $j\geq m$.
Combined with (4) this yields the desired lower bounds on
the partial degrees of the monomials $(q^{\si})^i(q^{\tau})^j$ occurring in 
$a^r_i(q^\si,q^\tau)$. An explicit expression for $a_r^r(q^{\si},q^{\tau})$ can
be given using Proposition \ref{th:intertwining}. The last statement
of the proposition follows then immediately.
\end{proof}
As a corollary we obtain the following crucial lemma.
\begin{lem}\label{th:limitlemma}
Let $-\infty<\tau<\infty$. The limits 
\begin{equation*}
\lim_{\si\to\infty} q^{r\si} e_r^{\si,\tau}, \quad 
\lim_{\si\to\infty} q^{2r\si}e_r^{\si,\si}
\quad (1\leq r\leq l)
\end{equation*}
exist in $\AU$. In other words, the coefficients of
$q^{r\si} e_r^{\si,\tau}$ respectively 
$q^{2r\si} e_r^{\si,\si}$ in the expansion with respect to the monomial basis 
of $\AU$ tend to finite values in the limit $\si\to\infty$.
\end{lem}
\begin{proof}
Fix $1\leq r\leq l$ and let  $-\infty < \sigma, \tau < \infty$. 
{}From Lemma \ref{th:limitlemma2} it is readily deduced that
\begin{equation*}
q^{r\sigma+r\tau}e_r =  
b^r_r(q^\sigma, q^\tau){\varphi_r^{\si,\tau}}_{|\TT}
+\cdots + b^r_0(q^\sigma, q^\tau){\varphi_0^{\si,\tau}}_{|\TT},
\end{equation*}
with $b_i^r$ ($0\leq i\leq r$) some polynomial in two variables 
and $b_r^r$ a non-zero constant polynomial (the important fact here
is that $b_i^r$ is a polynomial and not a Laurent polynomial).
Hence
\begin{equation}\label{est-polynomial}
e_r^{\si,\tau} = q^{-r\si-r\tau}\left 
(b_r^r(q^\si, q^\tau){\varphi_r^{\si,\tau}}
+ \cdots + b_0^r(q^\si, q^\tau){\varphi_0^{\si,\tau}}\right ).
\end{equation}
Since $\varphi_i^{\si,\tau} \to \varphi_i^{\infty,\tau}$ and 
$\varphi_i^{\si,\si}\to \varphi_i^{\infty,\infty}$ when 
$\si\to\infty$, the lemma follows.
\end{proof}
In view of Lemma \ref{th:limitlemma} we may set
for $1\leq r\leq l$ and $-\infty < \tau < \infty$,
\begin{equation}\label{e:einfdef}
\tilde{e}_r^{\infty,\tau} := \lim_{\si\to\infty} q^{r(\si+\tau -1)} 
(-1)^re_r^{\si,\tau},
\quad \tilde{e}_r^{\infty,\infty} := \lim_{\si\to\infty} q^{r(2\si -1)}
(-1)^re_r^{\si,\si}.
\end{equation}
It is clear from the definitions that $\tilde{e}_r^{\infty,\tau}\in \Hinft$ 
($1\leq r\leq l$, $-\infty < \tau \leq \infty$).
Observe that the elements $(-1)^re_r^{\si,\tau}$ are mapped onto $e_r(-x)\in \CC[x^{\pm
1}]^{\FSW}$ under the restriction mapping ${}_{|\TT}$.
 
Recall that the elements $\tilde{e}_r$ ($1\leq r\leq l$)
are algebraically independent generators of the algebra $\CC[x]^{\goS}$
(cf.\  \S \ref{section:polynomials}).
Again, we obtain an algebra isomorphism $P\mapsto \hat{P}$ of 
$\CC[x]^{\goS}$ onto $\CC[y]$, where
\begin{equation}\label{Sgeval}
\hat{P}(\tilde{e}_1(x),\ldots,\tilde{e}_l(x)):=P(x),\quad (P\in \CC[x]^{\goS}).
\end{equation}
\begin{thm}\label{th:limitBJ} Let $-\infty<\tau<\infty$.
The elements $\tilde{e}_r^{\infty,\tau}$ \textup{(}$1\leq r\leq l$\textup{)} 
mutually commute and
are algebraically independent generators of the algebra $\FSH^{\infty,\tau}$.
Any zonal $(\infty,\tau)$-spherical function 
$\varphi^{\infty,\tau}(\lambda)\in \FSH^{\infty,\tau}(\lambda)$ 
\textup{(}$\lambda\in P_K^+$\textup{)}
is equal to a non-zero scalar multiple of
\[
\hat{P}_{\lambda^{\natural}}^B(\tilde{e}_1^{\infty,\tau},\ldots,
\tilde{e}_l^{\infty,\tau};1,q^{2(n-2l)},1,q^{2\tau+2(n-2l)};q^2,q^2),
\]
where $P_{\mu}^B(.;a,b,c,d;q,t)$ is the multivariable big $q$-Jacobi polynomial
of degree $\mu$ \textup{(}cf.\ \S \ref{section:polynomials}\textup{)}.
\end{thm}
\begin{proof}
The elements $e_r^{\si,\tau}$ ($1\leq r\leq l$) mutually commute 
for $\si$ finite by Corollary \ref{commuterende}.  By the definition
of the elements $\tilde{e}_r^{\infty,\tau}$ \eqref{e:einfdef}, 
it follows that the $\tilde{e}_r^{\infty,\tau}$
($1\leq r\leq l$) also commute. Hence the element
$Q(\tilde{e}_1^{\infty,\tau},
\ldots,\tilde{e}_l^{\infty,\tau})\in \FSH^{\infty,\tau}$
for a polynomial $Q\in \CC[y]$ is well defined.

For $\si$ finite and $\lambda\in P_K^+$, 
a zonal $(\sigma,\tau)$-spherical function 
$\varphi^{\si,\tau}(\lambda)\in \Hst(\lambda)$ is given by
\begin{equation}\label{refpunt}
\varphi^{\si,\tau}(\lambda):=(-q^{\si+\tau-1})^{|\lambda^\natural|}
\hat{P}_{\lambda^\natural}^K(e_1^{\si,\tau},\ldots,
e_l^{\si,\tau};\underline{t}^{\si,\tau};q^2,q^2)
\end{equation}
with $P_{\mu}^K(.)$ the Koornwinder polynomial of degree $\mu$
(cf.\ Theorem \ref{th:AWspher}). 
Using the elementary properties of the Koornwinder polynomials given in Lemma
\ref{symmetryproperties},
this zonal $(\si,\tau)$-spherical function can be rewritten as
\[\varphi^{\si,\tau}(\lambda)=(s_{\ep})^{-|\lambda^\natural|}
\hat{P}_{\lambda^\natural}^K\bigl(s_{\ep}f_1^{\ep},\ldots, s_{\ep}^lf_l^{\ep};
\underline{t}_B(\ep);q^2,q^2\bigr),\]
where $s_\ep:= q/\ep (cd)^{\frac{1}{2}}$, 
$f_r^{\ep}:=(-s_{\ep})^{-r}e_r^{\si,\tau}$, $\ep:= q^{\si-(n-2l)}$, and
\[
\underline{t}_B(\ep)=\bigl(
\ep^{-1}(q^2c/d)^{\frac{1}{2}}, -\ep^{-1}(q^2d/c)^{\frac{1}{2}},
\ep a(q^2d/c)^{\frac{1}{2}},-\ep b(q^2c/d)^{\frac{1}{2}}\bigr),\]
with parameters $a,b,c,d$ given by
\begin{equation}\label{par}
a:=1, \; b:=q^{2(n-2l)}, \; 
c:=1, \; d:=q^{2\tau+2(n-2l)}.
\end{equation}
Observe that $s_{\ep}=q^{1-\si-\tau}$, hence by the definition \eqref{e:einfdef}
of $\tilde{e}_r^{\infty,\tau}$, $\lim_{\ep\downarrow 0}
f_r^{\ep}=\tilde{e}_r^{\infty,\tau}$ for all $r$. 
Combined with the limit transition 
from Koornwinder polynomials to multivariable big $q$-Jacobi polynomials 
\eqref{e:limitBJ-rewritten}, 
we get that $\varphi^{\infty,\tau}(\lambda):=
\lim_{\si\rightarrow\infty}\varphi^{\si,\tau}(\lambda)$ 
exists as limit in $\AU$,
and that
\begin{equation}\label{limitzonalB}
\varphi^{\infty,\tau}(\lambda)=
\hat{P}_{\lambda^{\natural}}^B(\tilde{e}_1^{\infty,\tau},\ldots,
\tilde{e}_l^{\infty,\tau};1,q^{2(n-2l)},1,q^{2\tau+2(n-2l)};q^2,q^2)
\end{equation}
with $P_{\mu}^B(.)$ the multivariable big $q$-Jacobi polynomial of degree $\mu$.
It is clear that $\varphi^{\infty,\tau}(\lambda)\in \FSH^{\infty,\tau}(\lambda)$,
but it may be zero since the algebraic independence 
of the elements $\tilde{e}_r^{\infty,\tau}$ ($r\in [1,l]$)
has not yet been established. 
To prove that $\varphi^{\infty,\tau}(\lambda)$ is non-zero, we compute the quadratic
norm $\|\varphi^{\si,\tau}(\lambda)\|^2$ with respect to the inner product
$\langle \varphi,\psi\rangle:=h(\psi^\ast\varphi)$, where $h$ is the normalized Haar
functional. Since all highest weights
$\lambda\in P_K^+$ are self-dual (i.e.\ $V(\lambda)$ is 
isomorphic to its dual representation), and since the two-sided coideal $\ksi$ is 
$\tau$-invariant, we have $\bigl(\varphi^{\si,\tau}(\lambda)\bigr)^\ast=
\varphi^{\si,\tau}(\lambda)$. Then it follows from the definition
\eqref{refpunt} of $\varphi^{\si,\tau}(\lambda)$, 
Corollary \ref{haarAW} and the definition \eqref{NX} of $N_K$, 
that
\begin{equation}\label{kwadraatnorm}
\|\varphi^{\si,\tau}(\lambda)\|^2=s_{\ep}^{-2|\lambda^\natural|}
N_K(\lambda^\natural;\underline{t}_B(\ep);q^2,q^2).
\end{equation}
The limit $\ep\downarrow 0$ (equivalently, $\si\to\infty$) 
in \eqref{kwadraatnorm}
can now be computed in the left-hand side and in the right-hand side 
(cf.\  Proposition \ref{normlimit}). It follows that
\[\|\varphi^{\infty,\tau}(\lambda)\|^2=N_B(\lambda^\natural;
1,q^{2(n-2l)},1,q^{2\tau+2(n-2l)};q^2,q^2).
\]
Since $N_B(\lambda^\natural)$ is strictly positive, it follows that
the quadratic norm 
$\|\varphi^{\infty,\tau}(\lambda)\|^2$ 
is non-zero, hence $\varphi^{\infty,\tau}(\lambda)\not=0$ for all $\lambda\in P_K^+$.
Hence the elements $\varphi^{\infty,\tau}(\lambda)\in\FSH^{\infty,\tau}(\lambda)$  
are zonal $(\infty,\tau)$-spherical functions for all $\lambda\in P_K^+$.

It remains to prove that the $\tilde{e}_r^{\infty,\tau}$ ($1\leq r\leq l$)
are algebraically independent. 
Consider the finite-dimensional subspaces
\[\FSH_m:=\bigoplus_{\lambda\in P_K^+: \lambda\leq m\varpi_l}
\FSH^{\infty,\tau}(\lambda),\quad (m\in {\mathbb{Z}}_+).\]
The dimension of the linear subspace $\FSH_m$ is equal to the number of 
positive integers $\underline{m}=(m_1,\ldots,m_l)\in {\mathbb{Z}}_+^{\times l}$ with
$|{\underline{m}}|:=\sum_im_i\leq m$, since $\varpi_r\leq\varpi_l$ for all $r\in
[0,l]$. For such a sequence of positive integers
$\underline{m}$, set $Q_{\underline{m}}(y):=
y_1^{m_1}\ldots y_l^{m_l}$. Since $\tilde{e}_r^{\infty,\tau}\in \oplus_{s=0}^r
\FSH^{\infty,\tau}(\varpi_s)$ and
\[
\FSH^{\infty,\tau}(\mu).\FSH^{\infty,\tau}(\mu')\subseteq
\bigoplus_{\nu\in P_K^+: \nu\leq\mu+\mu'}\FSH^{\infty,\tau}(\nu),\quad 
(\mu,\mu'\in P_{K}^+)\]
we have $Q_{\underline{m}}(\tilde{e}_1^{\infty,\tau},\ldots,\tilde{e}_l^{\infty,\tau})
\in\FSH_m$ for all $\underline{m}$ with $|\underline{m}|\leq m$. 
Hence the algebraic independence of $\tilde{e}_r^{\infty,\tau}$
($1\leq r\leq l$) will follow from the fact that the monomials
$Q_{\underline{m}}(\tilde{e}_1^{\infty,\tau},\ldots,e_l^{\infty,\tau})$ 
($|\underline{m}|\leq m$) span $\FSH_m$ for all $m\in {\mathbb{Z}}_+$.

Observe that $\FSH_m$ is spanned by the zonal spherical
functions $\varphi^{\infty,\tau}(\lambda)$ ($\lambda\leq m\omega_l$).
Since $P_{\mu}^B(x)$  ($\mu\in P_{\Sigma}^+$) is of the form 
$\tilde{m}_{\mu}(x)+ \sum_{\nu\in P_{\Sigma}^+: \nu<\mu}
c_{\nu}\tilde{m}_{\nu}(x)$ for certain constants $c_{\nu}$,
it follows from the explicit expression \eqref{limitzonalB} for 
$\varphi^{\infty,\tau}(\lambda)$ 
that each $\varphi^{\infty,\tau}(\lambda)$ with $\lambda\leq m\omega_l$
can be written as a linear combination of the monomials
$Q_{\underline{m}}(\tilde{e}_1^{\infty,\tau},\ldots,\tilde{e}_l^{\infty,\tau})$ 
($|\underline{m}|\leq m$). Hence,
the monomials
$Q_{\underline{m}}(\tilde{e}_1^{\infty,\tau},\ldots,\tilde{e}_l^{\infty,\tau})$ 
($|\underline{m}|\leq m$) span $\FSH_m$.
\end{proof}
\begin{thm}\label{th:limitLJ}
The elements $\tilde{e}_r^{\infty,\infty}$ 
\textup{(}$1\leq r\leq l$\textup{)} mutually commute and are algebraically
independent generators of the algebra $\FSH^{\infty,\infty}$.
Any zonal spherical function $\varphi^{\infty,\infty}(\lambda)$ 
\textup{(}$\lambda\in P^+_K$\textup{)} 
is equal to a non-zero scalar multiple of
\begin{equation*}
\hat{P}_{\lambda^\natural}^L
(\tilde{e}_1^{\infty,\infty},\ldots,\tilde{e}_l^{\infty,\infty};
q^{2(n-2l)}, 1;q^2,q^2)
\end{equation*}
where $P_{\mu}^L(.)$ is the multivariable little $q$-Jacobi polynomial of degree $\mu$
\textup{(}cf.\  \S \ref{section:polynomials}\textup{)}.
\end{thm}
\begin{proof}
By Theorem \ref{th:AWspher} and by symmetry properties of the Koornwinder
polynomials (cf.\  Lemma \ref{symmetryproperties}), a zonal
$(\si,\si)$-spherical function $\varphi^{\si,\si}(\lambda)\in
\FSH^{\si,\si}(\lambda)$ ($\lambda\in P_K^+$) is explicitly given
by
\[\varphi^{\si,\si}(\lambda):=(s_{\ep})^{-|\lambda^\natural|}
\hat{P}_{\lambda^\natural}^K(s_{\ep}f_1^\ep,\ldots,s_{\ep}^lf_l^\ep;
\underline{t}_L(\ep);q^2,q^2),\]
where $P_{\mu}^K(.)$ is the Koornwinder polynomial of degree $\mu$
(cf.\  \S \ref{section:polynomials}) and 
$\ep:= q^{2\si}$, $s_\ep:= q/\ep$, $f_r^\ep:= (-s_\ep)^{-r} e_r^{\si,\si}$, 
and $\underline{t}_L(\ep):=(\ep^{-1}q,-aq,\ep bq,-q)$ with $a:= q^{2(n-2l)}$ 
and $b:= 1$. Using \eqref{e:limitLJ-rewritten} together with the observation
that $\lim_{\ep\downarrow 0}f_r^\ep=\tilde{e}_r^{\infty,\infty}$, 
the proof is analogous to the proof of Theorem \ref{th:limitBJ}.  
\end{proof}
\begin{rem}
Using the limits $\tilde{e}_r^{\infty,\infty}=
\lim_{\tau\rightarrow\infty}\tilde{e}_r^{\infty,\tau}$ ($1\leq r\leq l$),
Theorem \ref{th:limitLJ} can also be proved 
by sending $\tau\to\infty$ in the results of Theorem \ref{th:limitBJ}.
On the level of multivariable orthogonal polynomials this limit corresponds
to the limit from multivariable big $q$-Jacobi polynomials 
to multivariable little $q$-Jacobi polynomials proved in
\cite[Thm. 5.1(3)]{stok-koor:limit}.
\end{rem}
\begin{rem}
As a corollary of Theorem \ref{th:limitBJ} and Theorem \ref{th:limitLJ},
the restricted Haar functional $h: \Hst\rightarrow \CC$
for $\si=\infty$, $-\infty<\tau<\infty$ 
respectively for $\sigma=\tau=\infty$  can be expressed
in terms of the orthogonality 
measure of the multivariable big respectively 
little $q$-Jacobi polynomials (cf.\ Corollary \ref{haarAW}).
\end{rem}
\begin{rem}
In the last two sections 
we have interpreted 
the Koornwinder polynomial resp.\ the multivariable big and the multivariable 
little $q$-Jacobi polynomial
\begin{equation}
\begin{split}
&P_{\mu}^K(.\,;-q^{\si+\tau+1},-q^{-\si-\tau+1},q^{\si-\tau+1},
q^{-\si+\tau+2(n-2l)+1};q^2,q^2),\\ 
&P_{\mu}^B(.\,;1,q^{2(n-2l)},1,q^{2\tau+2(n-2l)};q^2,q^2),\\
&P_{\mu}^L(.\,;q^{2(n-2l)},1;q^2,q^2)
\end{split}
\end{equation}
($-\infty<\si,\tau<\infty$, $\mu\in P_K^+$)
as a zonal spherical function 
on some quantum analogue of the complex Grassmannian.
For each of these polynomials, the classical limit $q\uparrow 1$ can be computed
using results from \cite[\S 4]{dj:polynomial} for the Koornwinder polynomials
and using \cite[Thm. 5.1(5)\&(6)]{stok-koor:limit} for the multivariable big and
little $q$-Jacobi polynomials.
The limits can be written in terms of the generalized Jacobi polynomial 
$P_{\mu}^J(.\,;n-2l,0;1)$ or, equivalently, in terms of 
the $BC$ type Heckman-Opdam polynomial
$P_{\mu}^{HO}(.\,;n-2l,1,1/2)$.
This agrees nicely with the classical interpretation of 
the Heckman-Opdam polynomial $P_{\mu}^{HO}(.\,;n-2l,1,1/2)$ 
as a zonal spherical function on the complex Grassmannian $U/K$ (see \S 2).
\end{rem}

\section*{Acknowledgements}
The authors would like to thank Prof.\ Tom H. Koornwinder and  
Prof.\ Masatoshi Noumi for interesting and informative discussions 
about the topic of this paper.
%%%%%%%%%%%%%%%%%%%%%%%%%%%%%%%%%%%%%%%%%%%%%%%%%%%%%%%%%%%%%%%%%%%%%%%%%%%%%%%%%%
%%
%%
%%   REFERENCES
%%
%%
%%%%%%%%%%%%%%%%%%%%%%%%%%%%%%%%%%%%%%%%%%%%%%%%%%%%%%%%%%%%%%%%%%%%%%%%%%%%%%%%%%
\bibliographystyle{plain}

\end{document}